# Primes in quadratic fields

Theodorus J. Dekker [*])

## Abstract

This paper presents algorithms for calculating the quadratic character and the norms of prime ideals in the ring of integers of any quadratic field. The norms of prime ideals are obtained by means of a sieve algorithm using the quadratic character for the field considered.

A quadratic field, and its ring of integers, can be represented naturally in a plane. Using such a representation, the prime numbers - which generate the principal prime ideals in the ring - are displayed in a given bounded region of the plane.

## 1. Introduction

This paper presents an algorithm to calculate the quadratic character for any quadratic field, and an algorithm to calculate the norms of prime ideals in its ring of integers up to a given maximum. The latter algorithm is used to produce pictures displaying prime numbers in a given bounded region of the ring.

A quadratic field, and its ring of integers, can be displayed in a plane in a natural way. The presented algorithms produce a complete picture of its prime numbers within a given region.

If the ring of integers does not have the unique-factorization property, then the picture produced displays only the prime numbers, but not those irreducible numbers in the ring which are not primes. Moreover, the prime structure of the non-principal ideals is not displayed except in some pictures for rings of class number 2 or 3.

Section 2 gives the theory needed on algebraic number fields, in particular quadratic fields, prime numbers and ideals, and quadratic characters. The theory presented in this section is not new. For details and proofs see, for instance, [Lejeune Dirichlet / Dedekind 1893], [Borewicz & Safarevic 1966], [Cohn 1978], [Ireland & Rosen 1982], [Cohen 1993].

In section 3 we present our algorithms and in appendices we present several pictures displaying prime numbers and, for some rings also non-principal prime ideals.

In a previous paper, [Dekker 1994], we presented a similar algorithm to calculate (norms of) prime numbers and to produce a picture displaying prime numbers for a certain field. In that paper the quadratic characters needed were calculated by hand for some specific example fields.

---

*) address: Van Uytrechtlaan 25, NL-1901 JK Castricum; e-mail: T.J.Dekker@uva.nl .

## 2. Theory

2.1. Algebraic number fields

An 'algebraic number' is a zero of a rational polynomial, i.e. a polynomial whose coefficients are elements of field $\mathbb{Q}$ of rational numbers. An 'algebraic number field' is an extension field of finite degree over $\mathbb{Q}$, which is obtained from $\mathbb{Q}$ by adjoining an algebraic number.

An 'algebraic integer' is a zero of an integral monic polynomial, i.e. a polynomial whose coefficients are elements of ring $\mathbb{Z}$ of integers and whose leading coefficient is 1. The algebraic integers in an algebraic number field form an 'integral domain', i.e. a commutative ring with identity and without zero-divisors; this domain is called the 'ring of integers' or 'domain of integers' of the field, in the sequel denoted by R.

Elements of $\mathbb{Z}$ are often called 'rational integers' to avoid confusion with algebraic integers. Lower case latin letters denote rational numbers or rational integers, lower case Greek letters denote algebraic numbers or algebraic integers, and upper case letters denote rings or ideals.

Units, prime numbers and irreducibles

A 'unit' of an algebraic number field (or of its ring of integers) is an integer of the field whose inverse is also an integer of the field.

Two integers $\alpha$ and $\beta$ of an algebraic number field are called 'associates' of each other whenever $\beta = \varepsilon \alpha$, where $\varepsilon$ is a a unit of the field.

An integer of an algebraic number field is called 'composite', if it is the product of two non-unit integers of the field. An integer of an algebraic number field is called 'irreducible', if it is not a unit and not composite. Thus, within the ring of integers of a field, an irreducible number is divisible only by the element itself and its associates and by the units of the field.

Note that divisibility is always understood within the ring of integers of the field considered. Thus, for integers $\alpha$ and $\beta$ we define: $\beta$ 'is divisible by' $\alpha$, (or, equivalently, $\alpha$ 'divides' $\beta$, notation $\alpha \mid \beta$), if there is an integer $\gamma$, such that $\alpha\gamma = \beta$.

A non-zero integer $\pi$ of an algebraic number field is called a 'prime' (of the field or of its ring of integers), if it is not a unit, and it satisfies the property that whenever $\pi \mid \alpha\beta$ for any integers $\alpha$ and $\beta$, then $\pi \mid \alpha$ or $\pi \mid \beta$.

Unique factorization

A ring (integral domain) is said to have the 'unique factorization property', or to be a 'unique factorization ring' ('unique factorization domain', abbreviated 'UFD'), if the factorization of any integer into irreducible factors is unique apart from the order of the factors and ambiguities between associated factors, i.e. the irreducible factors in any two factorizations of a given integer can be ordered such that corresponding factors are associates.

The distinction between primes and irreducibles is needed, because in a ring which is not an UFD these notions are not equivalent. Every prime is irreducible, but the converse is not true in general. In a UFD, however, every irreducible number is prime.

The domain of integers of an algebraic number field mostly is not a UFD. Instead of numbers as products of irreducibles, we consider ideals as products of prime ideals.

Ideals

An 'ideal' A in the ring R of integers of an algebraic number field is kernel of a homomorphism from R to another ring; it is a subset of R which contains, with any elements $\alpha$ and $\beta$ of A, also the difference $\alpha - \beta$ and the product $\alpha\rho$, where $\rho$ is any element of R.

An ideal generated by one element α of R, notation (α), is called a 'principal ideal'. In particular we have the zero ideal (0) containing only the number zero, and the unit ideal (1), which is equal to R.

If all ideals in R are principal ideals, then R is called a 'principal ideal domain'. A principal ideal domain necessarily is a unique-factorization domain.

An ideal A of R induces subdivision of R into residue classes: two elements of R belong to the same residue class if and only if their difference belongs to A. The residue classes form a ring isomorphic to the homomorphic image R/A of R. Moreover, this ring satisfies:

2.1.1. Theorem. If R is the ring of integers of an algebraic number field and A is a non-zero ideal in R, then R/A is finite.

This allows us to define the (absolute) 'norm' of A, notation N(A), as follows:

If A is not the zero ideal, then N(A) is equal to the size of R/A, otherwise N(A) equals 0.

The product of two ideals A and B is the ideal generated by all products αβ, where α is in A and β in B. Ideal A is said to be 'divisor of' B if there is an ideal C in R such that AC = B, or, equivalently, if A contains B.

2.1.2. Theorem. The norm of a product of two ideals equals the product of the norms.

A non-zero and non-unit ideal P in R is called 'prime ideal' if the ring R/P of its residue classes is an integral domain. Equivalently, a non-zero and non-unit ideal P is a prime ideal if it satisfies the property that whenever P | AB for any ideals A and B, then P | A or P | B. This property is similar to the definition of prime number given above. In fact, a principal ideal generated by a prime number is a prime ideal, but a principal ideal generated by an irreducible which is not prime, is not a prime ideal.

The ring R/P, being a finite integral domain, necessarily is a field; therefore its order is a power of a prime number. Hence:

2.1.3. Theorem. The norm of a prime ideal is a power of a prime number.

The ring of integers of an algebraic number field is a 'Dedekind ring', i.e. satisfies:

2.1.4. Theorem (Dedekind, 1893). Every non-zero ideal in the ring of integers of an algebraic number field can be written uniquely (up to order of the factors) as a product of prime ideals.

2.2. Quadratic fields

A 'quadratic field' is obtained from $\mathbb{Q}$ by adjoining a quadratic algebraic number, for which we can simply take √r, where r is a non-square (and square-free) rational integer. This field, denoted by $\mathbb{Q}$(√r), consists of the numbers ζ:= x+y√r, where x and y are rational. Number r is called the (square-free) 'radicand' of the field.

The algebraic integers of field $\mathbb{Q}$(√r) form a ring R = $\mathbb{Z}$[τ], obtained from $\mathbb{Z}$ by adjoining an appropriate algebraic integer τ of $\mathbb{Q}$(√r), i.e. $\mathbb{Z}$[τ] consists of the numbers x + τ y, where x and y are rational integers. The value of τ, depending on r, is given by

$$\tau = \sqrt{r}, \quad \text{when } r \equiv 2 \text{ or } r \equiv 3 \pmod{4},$$
$$\tau = (1 + \sqrt{r})/2, \quad \text{when } r \equiv 1 \pmod{4}.$$

The case r ≡ 0 (mod 4) is excluded, since r is assumed square-free.

The 'discriminant' of an irrational element ζ of $\mathbb{Q}$(√r) is the discriminant of the irreducible rational (monic) polynomial of degree 2, having zero ζ. It is equal to r times the square of a rational number.

If ζ is an algebraic integer of the field, then its discriminant equals r times the square of a rational integer. The 'discriminant' of field $\mathbb{Q}(\sqrt{r})$ is defined as the discriminant of minimal magnitude of its irrational integers, i.e. it equals the discriminant of τ. So, this discriminant, denoted by d, has the following value:

$$d = r, \quad \text{when } r \equiv 1 \pmod{4};$$
$$d = 4r, \quad \text{when } r \equiv 2 \text{ or } r \equiv 3 \pmod{4}.$$

The (absolute) 'norm' N(ζ) of ζ is defined as the absolute value of the product of ζ and its conjugate, i.e.

$$N(\zeta) := |x^2 - r y^2|.$$

In the ring of algebraic integers $\mathbb{Z}[\tau]$, the norm of ζ = x + τ y is similarly as follows:

$$N(\zeta) = |x^2 - r y^2|, \quad \text{when } r \equiv 2 \text{ or } r \equiv 3 \pmod{4},$$
$$N(\zeta) = |x^2 + x y - c y^2|, \ c := (r-1)/4, \quad \text{when } r \equiv 1 \pmod{4}.$$

This definition corresponds to the definition of norm of ideals given above as follows.

<u>2.2.1. Theorem.</u> Let A = (ζ) be the principal ideal generated by ζ in $\mathbb{Z}[\tau]$. Then N(ζ) = N(A).

Norms have the following properties.

<u>2.2.2. Theorem.</u> The norm of a product of two elements equals the product of the norms. The norm of an algebraic integer of the field is an integer.

Furthermore, we have

<u>2.2.3. Theorem.</u> (cf. Theorem 2.1.3). The norm of a prime number, as well as the norm of a prime ideal, is either a prime number or the square of a prime number.

<u>2.2.4. Theorem.</u> If ζ is an integer of a quadratic field and N(ζ) is a prime in $\mathbb{Z}$, then ζ is prime.

<u>2.2.5. Theorem.</u> Every integer of a quadratic field, not zero or a unit, can be factored into a finite product of irreducible integers of the field.

Quadratic fields and unique factorization

For several quadratic fields, the ring of integers is not a unique-factorization domain. In that case, the ring contains irreducible numbers which are not primes.

There exist only a finite number of complex quadratic fields (r < 0) whose ring of integers is a unique-factorization domain, namely the fields whose radicand values are:

$$-1, -2, -3, -7, -11, -19, -43, -67, -163.$$

The theorem that no other complex quadratic fields satisfy has been proved independently by [Baker 1966] and [Stark 1967], the latter correcting an error in the proof by [Heegner 1952].

For real quadratic fields (r > 0) it is unknown if there are infinitely many fields having unique factorization. Up to r < 100 unique factorization holds exactly for the following values:

r ≡ 1 (mod 4):       5, 13, 17, 21, 29, 33, 37, 41, 53, 57, 61, 69, 73, 77, 89, 93, 97;

r ≡ 2 (mod 4):       2, 6, 14, 22, 38, 46, 62, 86, 94;

r ≡ 3 (mod 4):       3, 7, 11, 19, 23, 31, 43, 47, 59, 67, 71, 83.

For surveys of these and related results see, for instance, [Borewicz & Safarevic 1966, Tabelle p. 454] and [Cohen 1993].

[Behrbohm & Rédei 1936] showed (see also [Van der Linden 1984]) that a real quadratic field can have the unique-factorization property only if the radicand r has at most two different prime factors; either r is prime or r = pq, where p ≡ 3 (mod 4) and q = 2 or q ≡ 3 (mod 4).

The fields $\mathbb{Q}(\sqrt{r})$, for $0 < r < 100$, satisfying this criterion, all have the unique-factorization property with only one exception, namely r = 79.

A consequence of this theorem is that r = 2 is the only value of r ≡ 2 (mod 8) for which the field $\mathbb{Q}(\sqrt{r})$ has the unique-factorization property.

### 2.3. Prime ideals in quadratic fields

Let $R = \mathbb{Z}[\tau]$ be the ring of integers of quadratic field $\mathbb{Q}(\sqrt{r}) = \mathbb{Q}(\sqrt{d})$, where r is its square-free radicand and d its discriminant. To determine if an ideal in R is a prime ideal or not, depends only on its norm. This follows from the following theorems.

#### 2.3.1. Theorem. For any prime ideal P in R, a unique natural prime p exists such that (p) | P.

For the next theorems we need the concept of 'quadratic residue' defined as follows.
Let a be a non-zero rational integer and p a natural prime not divisor of a. Then a is called 'quadratic residue' modulo p if there is a natural number n such that $n^2 \equiv a$ (mod p); otherwise, a is called 'quadratic non-residue' modulo p.

The set of non-zero residues modulo p consists of the p-1 elements 1, ..., p-1. If p is odd, then this set contains (p-1)/2 quadratic residues and the same number of quadratic non-residues.

#### 2.3.2. Theorem. Let p be a natural prime.
1. The ideal (p) in R is either a product of two prime ideals in R, both having norm equal to p, or is itself a prime ideal in R, having norm equal to $p^2$; applying this to all natural primes, all prime ideals in R are obtained;
2. If p (odd or even), is divisor of d, then (p) is product of two equal prime ideals in R;
3. If p is odd and not divisor of d (or r), then (p) is product of two distinct prime ideals if d is a quadratic residue modulo p; otherwise, (p) itself is a prime ideal in R;
4. If 2 is not divisor of d, hence d = r ≡ 1 (mod 4), then (2) is product of two distinct prime ideals in R if d ≡ 1 (mod 8); otherwise, d ≡ 5 (mod 8) and (2) itself is prime ideal in R.

A simpler equivalent formulation of this theorem is obtained by means of the 'Legendre-Jacobi-Kronecker' symbol, notation (a/b), where a and b are rational integers.

#### 2.3.3. Definition of Legendre-Jacobi-Kronecker symbol [Cohen 1993, section 1.4.2]:
For any rational integers a and b, symbol (a/b) is defined as follows:

1.   Legendre symbol for odd natural p:

   (a/p) = +1 if a is quadratic residue modulo p,

   (a/p) = 0 if p divides a,

   (a/p) = −1 if a is quadratic non-residue modulo p.

2.   Extension by Kronecker: for b = 2:

   (a/2) = +1 if a ≡ 1 or 7 (modulo 8),

   (a/2) = 0 if a is even,

   (a/2) = −1 if a ≡ 3 or 5 (modulo 8).

3. Jacobi's product rule for $b > 0$, $b = \Pi\, p$ is product of its – odd or even – prime factors p:

$$(a/b) = \Pi\, (a/p);$$

this includes the case $b = 1$, for which case $(a/b) = 1$.

4. Extension for $b \leq 0$,

$$(a/b) = (a/-b) \text{ if } b < 0 \text{ and } a \geq 0,$$

$$(a/b) = -(a/-b) \text{ if } b < 0 \text{ and } a < 0,$$

$$(a/0) = 1 \text{ if } a = \pm 1, \text{ otherwise } 0.$$

Remarks

1. In this paper we need the symbol $(a/b)$ only for $a = d$, where d is discriminant of a quadratic field.

2. If, for odd positive b, all factors $(a/p)$ of $(a/b)$ have the value +1, then a is quadratic residue modulo b, according to the Chinese remainder theorem. Thus, $(a/b) = 1$ is necessary, but certainly not sufficient, to ensure that a is quadratic residue modulo b.

2.3.4. Theorem of quadratic reciprocity (C.F. Gauss, 1796)

Let a and b be any odd natural numbers having no common prime factor. Then the Legendre-Jacobi-Kronecker symbol satisfies

$$(a/b) = -(b/a) \text{ if both a and b are congruent 3 modulo 4,}$$

$$(a/b) = (b/a) \text{ otherwise.}$$

This theorem enables us, in particular, to transform the property of quadratic residue of d (or r) modulo odd prime p to a property of odd primes p modulo |d|.

2.3.5. Definition of Kronecker's quadratic character (cfr. [Cohn 1978])
Consider a quadratic field $\mathbb{Q}(\sqrt{d})$, where d is the discriminant of the field. The 'quadratic character' of this field is the function $\chi(x) = \chi_d(x)$, defined, for any rational integer x, by:

$$\chi_d(x) := (d/x).$$

2.3.6. Theorem.
The quadratic character $\chi(x) = \chi_d(x)$ of a field $\mathbb{Q}(\sqrt{d})$ has the following properties:

1) it is multiplicative:   $\chi(xy) = \chi(x)\chi(y)$;
2) it is periodic modulo |d|:   if $x \equiv y \pmod{|d|}$ then $\chi(x) = \chi(y)$;
3) if $d > 0$ then it is symmetric:   $\chi(-x) = \chi(x)$,
   if $d < 0$ then it is anti-symmetric:   $\chi(-x) = -\chi(x)$.

Using only parts 1 and 2 of 2.3.3, theorem 2.3.2 is equivalent to the following simpler:

2.3.7. Theorem. Let $R = \mathbb{Z}[\tau]$ be the ring of integers of quadratic field $\mathbb{Q}(\sqrt{d})$ with discriminant d, and p a natural prime, odd or even. Then the following cases yield all prime ideals in R:
1. if $= \chi_d(p) = 1$, then p 'splits', i.e. (p) is product of two distinct prime ideals P, P' in R:
   $$(p) = P\, P',\ P \neq P',\ N(P) = N(P') = p;$$
2. if $= \chi_d(p) = 0$, then p 'ramifies', i.e. (p) is product of two equal prime ideals P in R:
   $$(p) = P\, P,\ N(P) = p;$$
3. if $= \chi_d(p) = -1$, then p is 'inert', i.e. (p) itself is a prime ideal in R:
   $$(p) = P,\ N(P) = p^2.$$

## 3.   Algorithms

Using the above theory, we now present our algorithms to calculate quadratic characters and norms of prime ideals which enable us to produce pictures of prime numbers and ideals in quadratic rings.

3.1. Algorithm to calculate quadratic character

To calculate the quadratic character of field $\mathbb{Q}(\sqrt{r})$, firstly the given number, r, is reduced to a square-free number by removing squares. Moreover, the prime factors of the resulting square-free radicand are obtained and the discriminant, d, is determined as follows.

For odd prime p, let p' = p if p $\equiv$ 1 (mod 4), otherwise p' = -p. Then we have for square-free r:

If r $\equiv$ 1 (mod 4) then d = r = $\Pi$ p', the product being taken for all odd prime factors p of r,

otherwise d = 4r = e $\Pi$ p', where e = -4 for r $\equiv$ 3 (mod 4) and e = +8 or -8 for r $\equiv$ 2 or 6 (mod 8).

Subsequently, the quadratic character $\chi_d(x)$ is calculated for one period $0 \le x < |d|$ as follows: for each odd prime factor p of d, the character $\chi_{p'}(x)$ is set to 0 for 0, to +1 for the non-zero squares modulo p and to -1 for the remaining numbers modulo p.

The total quadratic character modulo |d| is obtained, using the multiplicative property (see theorem 2.3.6), by multiplying the characters for the odd prime factors of d and - when d $\equiv$ 0 modulo 4 - of the remaining factor e = -4 or $\pm$ 8 defined above. These characters for the corresponding cases are (cf. definition 2.3.3 and theorem 2.3.7):

   if r $\equiv$ 3 (mod 4):     e = -4,        $\chi_e(x) = (0, +1, 0, -1)$,

   if r $\equiv$ 2 (mod 8):     e = +8,        $\chi_e(x) = (0, +1, 0, -1, 0, -1, 0, +1)$,

   if r $\equiv$ 6 (mod 8):     e = -8,        $\chi_e(x) = (0, +1, 0, +1, 0, -1, 0, -1)$.

The total number of operations for calculating the quadratic character is of the order of |d|.

3.2. Algorithm to determine norms of prime ideals

Whether an ideal in the ring of integers of a quadratic field is prime or not, depends only on its norm. This follows from theorems 2.3.1 and 2.3.2 (or 2.3.7). We propose, therefore, an algorithm to calculate a set of norms of prime ideals in the ring of integers, R, of a given quadratic field. From this, one easily obtains the corresponding set of prime numbers of the field.

For given discriminant d (or radicand r) of quadratic field $\mathbb{Q}(\sqrt{d})$ and a certain maximum value 'max', our algorithm determines the norms of all prime ideals of the field not larger than max.

First an appropriate set S = S(d, max) of natural numbers is formed, to be used as starting set for a sieve algorithm. S must at least contain the norms not larger than max of all prime ideals in R, namely the non-inert prime numbers and the squares of the inert primes. S may, to begin with, contain also norms of non-prime ideals, but must not contain any inert prime number, whose square is norm of a prime ideal (see theorem 2.3.7).

Subsequently, a sieve algorithm, defined below, is applied to S in order to remove the norms of all non-prime ideals. Thus all norms of prime ideals in R up to the maximum value max are obtained.

Sieve algorithm

We define an appropriate starting set S by means of the quadratic character defined above (2.3.5), which, according to theorem 2.3.6, is periodic modulo the discriminant of the field considered.

Let S = S(d, max) be the set of natural numbers n, $2 \leq n \leq$ max, satisfying either n is (prime) divisor of d, or otherwise $\chi_d(n) = 1$.

The 'sieve algorithm' applied to S yields a set consisting of those elements of S which are not product of two other elements of S. This is achieved by removing all elements which are products of elements of S, as follows.

Starting from set T := S, the algorithm proceeds in successive steps, where in each step certain elements are removed from T.

In each step, let t be the smallest element of T such that $\chi_d(t) = 1$ and t has not yet been treated in previous steps. Then we remove those elements of T which are product of t and some element of S. In fact, the products are removed starting from $t^2$, since the smaller products have already been removed in previous steps.

The algorithm is completed when the next element of T to be treated is larger than the square-root of max, and set T is then delivered as the result of the algorithm.

3.2.1. Theorem

Let $\mathbb{Q}(\sqrt{d})$ be the quadratic field having discriminant d, and R be its ring of integers.

Let S = S(d, max) be the set defined above.

Then the sieve algorithm, defined above, applied to S(d, max) yields a set T containing:

1. all prime numbers p satisfying $\chi_d(p) = 0$ or 1;

2. all squares of prime numbers satifying $\chi_d(p) = -1$;

3. all products pq of distinct primes numbers p and q satifying $\chi_d(p) = \chi_d(q) = -1$.

Hence, T contains the norms not larger than max of all prime ideals in R and no norms of non-prime ideals in R.

Proof. This follows from theorems of section 2.3.    []

Remarks

1. The third subset of T consists of numbers which are not norm of any ideal in R; therefore these numbers are superfluous but harmless.
2. For efficiency reasons, we apply the sieve process to odd numbers only. Thus, we start with set S containing, besides the prime divisors of d, the odd numbers satisfying $\chi_d(n) = 1$ and, when $d \equiv 1 \pmod 4$, also the unique even prime norm 2 or 4. In the latter case, the quadratic character of odd numbers is periodic modulo 2|d|.
3. The number of operations required for the sieve algorithm is of the order $\rho$ max $\Sigma_{p \leq \sqrt{max}} 1/p$, where $\rho$ is the density of (odd) numbers of quadratic character $\chi_d(n) = +1$ within a period |d| (or 2 |d|). This is asymptotically equal to $\rho$ max log log max and is a fraction $\rho$ of the number of operations required for the sieve of Eratosthenes to calculate natural primes. The density $\rho$ is equal to $\phi(|d|) / 2|d|$, where $\phi$ denotes Euler's function, i.e. the number of natural numbers smaller than and relative prime to |d|.

Programs and pictures

We have implemented our algorithms in C (sub-)programs to calculate quadratic characters, to calculate prime norms and to generate pictures of primes. We used the Code Warrior C-compiler of Metroworks Corporation including its drawing facility to produce the pictures.

For details and for more and higher-resolution pictures, see our website:

    http://www.science.uva.nl/~dirk .

**Appendices: pictures of primes**

In appendices we present pictures of prime numbers – and in the last four appendices also of non-principal prime ideals - for several complex and real quadratic fields.

At the top, each picture mentions the field, $\mathbb{Q}(\sqrt{r})$, and displays its quadratic character (as far as space allows). In the pictures, rational integers are placed on the x-axis and numbers of the form $\sqrt{r}$ times rational integers on the y-axis. When $d \equiv 0$ modulo 4, we use a square grid, otherwise a staggered grid, where the grid points form roughly equilateral triangles.

For quadratic non-UFD, the pictures display the prime numbers, which generate the principal prime ideals, but not those irreducible numbers which are not prime.
Moreover, only in appendices 5-8, the non-principal prime ideals are displayed as follows. The non-principal ideals are obtained by dividing principal ideals by a certain non-principal prime ideal, I, or its conjugate, where I is defined by (cf. section 2.2)

    $I := [\text{norm}, \zeta]$,    $\zeta := \text{shift} + \tau = \text{shift} + (d \bmod 4 + \sqrt{d}) / 2$,

i.e. I is generated by 'norm' being its norm, and the integer $\zeta$ of $Q(\sqrt{r})$.
In the picture, the non-principal prime ideals then are represented by those numbers whose norm is equal to a prime norm times the norm of I.

**Acknowledgments**

The author is grateful to Dr. Herman te Riele for reading the manuscript and for his valuable comments.

## Appendix 1:        Pictures of prime numbers for complex UFD

The pictures show the quadratic character and a picture of <span style="color:blue">prime numbers</span> and <span style="color:red">units</span> for the complex quadratic fields whose domain of integers is a unique-factorization domain, namely

the fields of discriminant congruent 0 modulo 4:

   $Q(\sqrt{-1}), Q(\sqrt{-2})$

and the fields of discriminant congruent 1 modulo 4:

   $Q(\sqrt{-3}), Q(\sqrt{-7}), Q(\sqrt{-11}), Q(\sqrt{-19}), Q(\sqrt{-43}), Q(\sqrt{-67}), Q(\sqrt{-163})$.

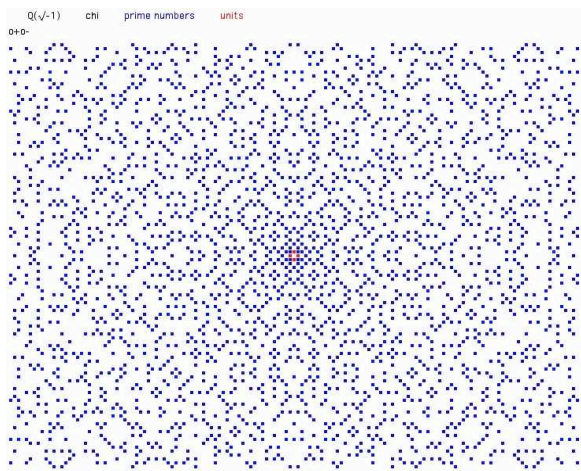
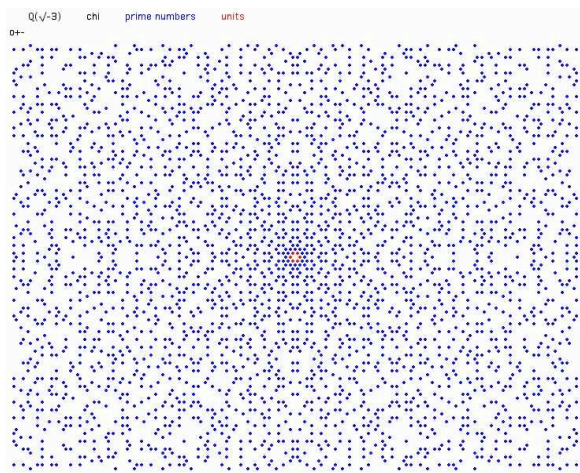
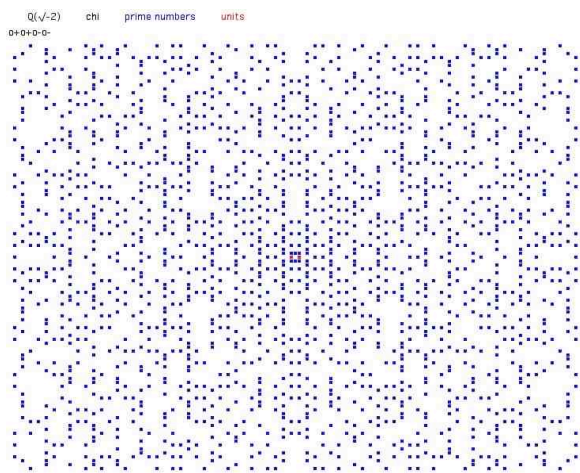
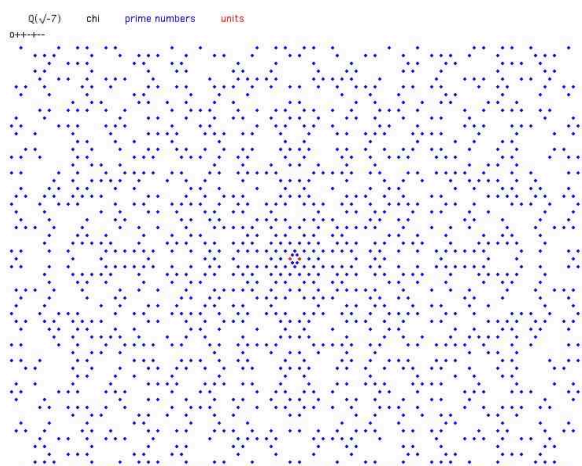

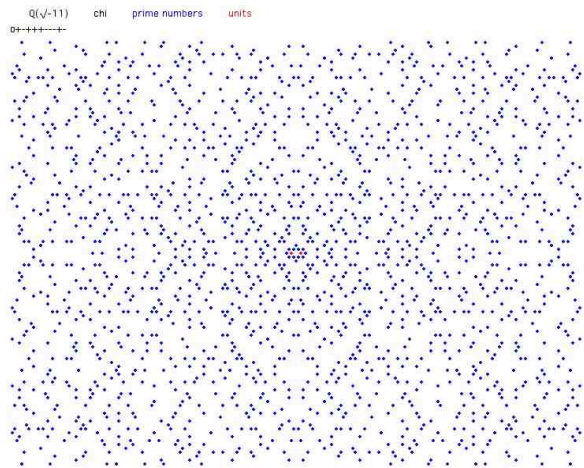
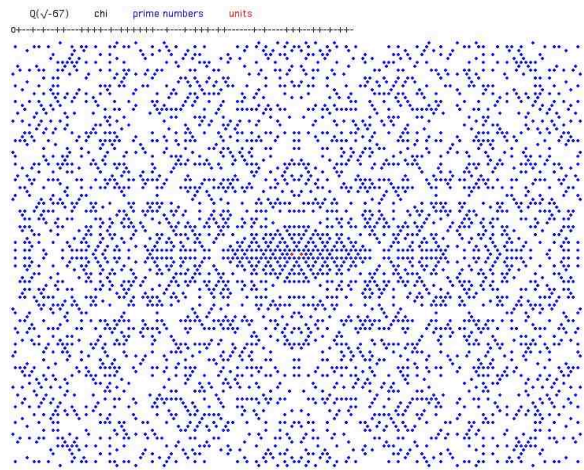
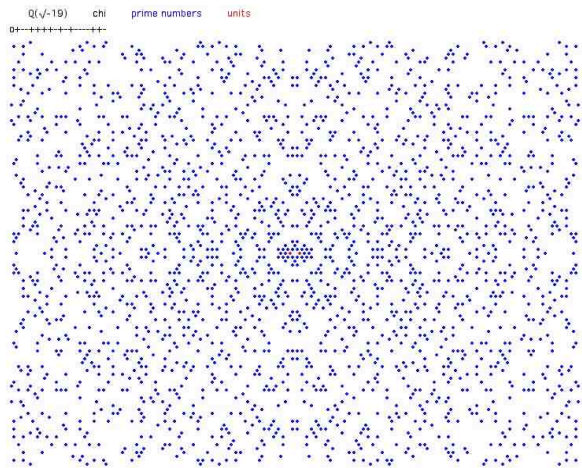
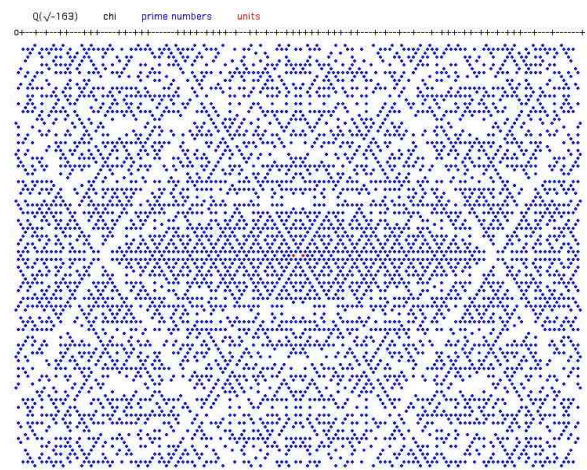
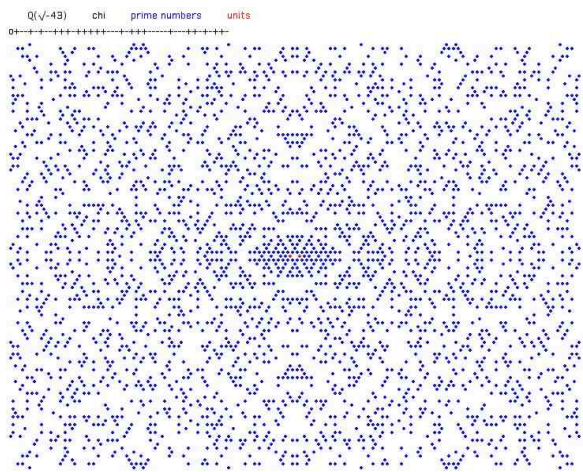

Appendix 2:            Pictures of prime numbers for real UFD

The pictures show the quadratic character and a picture of prime numbers and units for some real quadratic fields whose domain of integers is a unique-factorization domain, namely, for radicands < 32,

the fields of discriminant congruent 0 modulo 4:

Q(√2), Q(√3), Q(√6), Q(√7), Q(√11), Q(√14), Q(√19), (√22), Q(√23), Q(√31)

and the fields of discriminant congruent 1 modulo 4:

Q(√5), Q(√13), Q(√17), Q(√21), Q(√29).

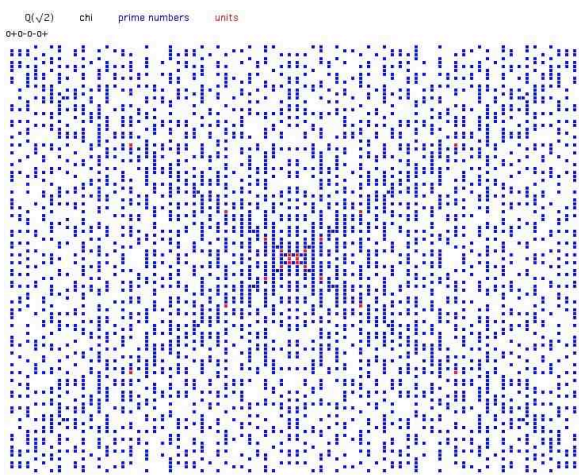
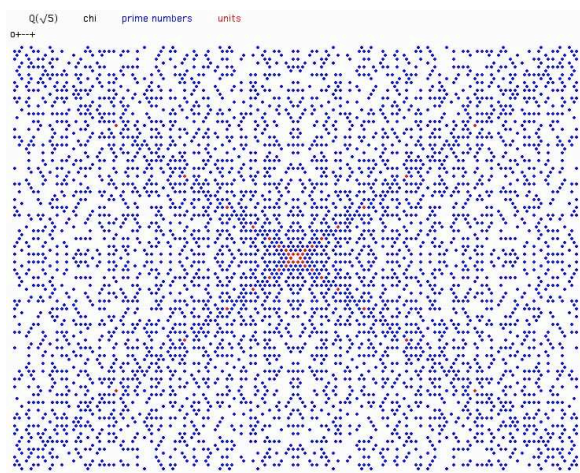
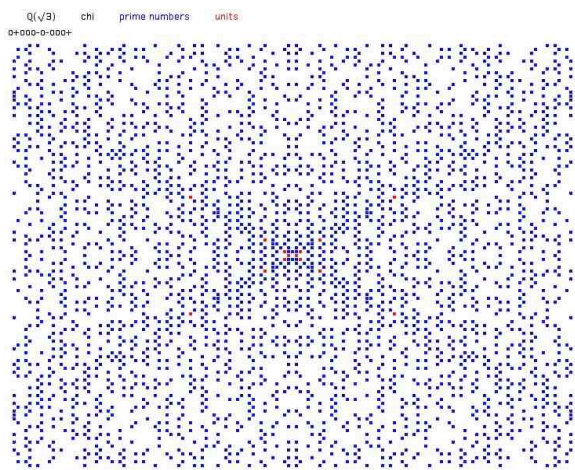
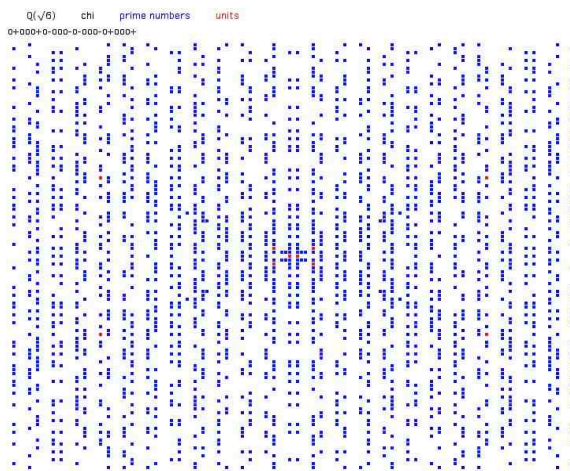

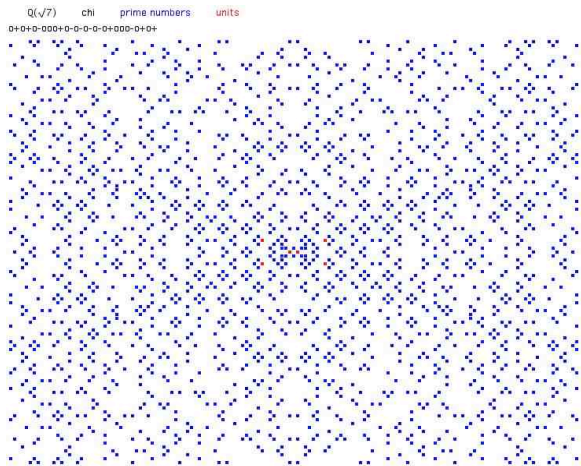
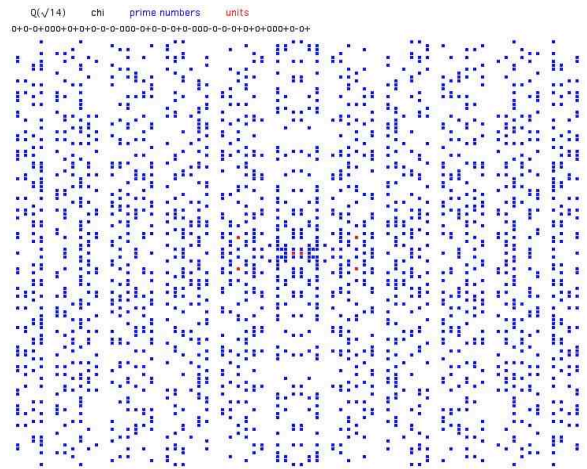
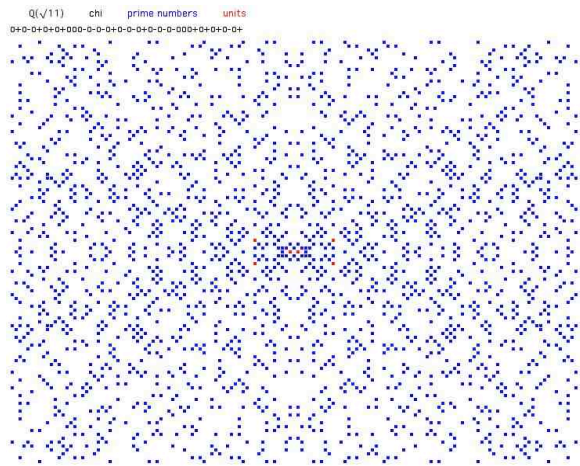
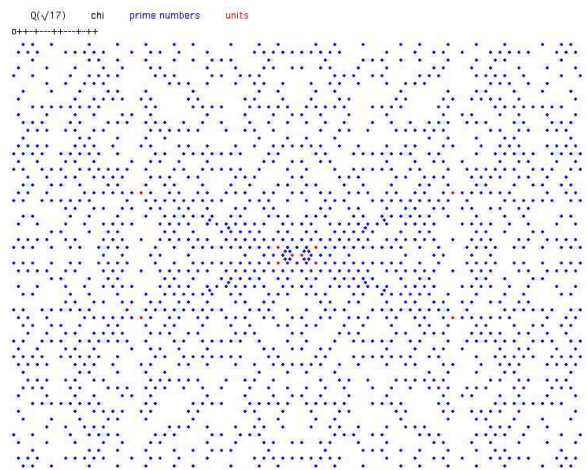
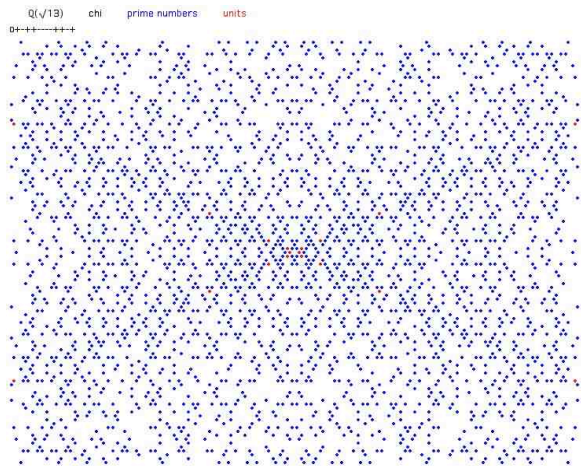
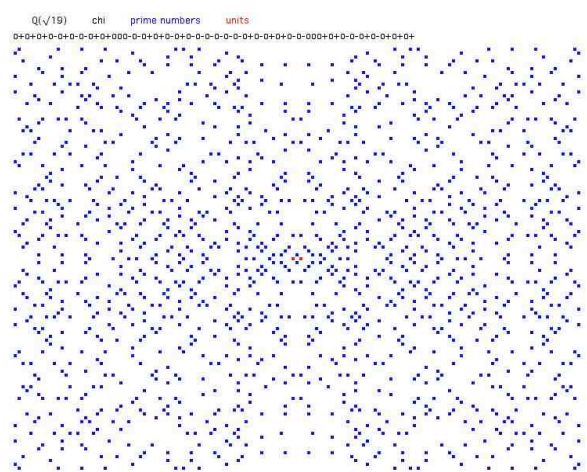

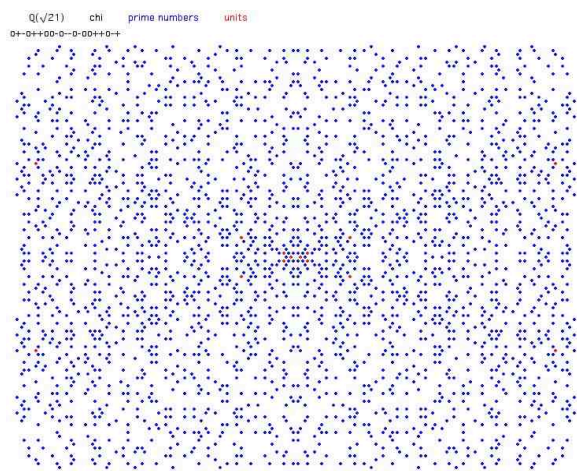
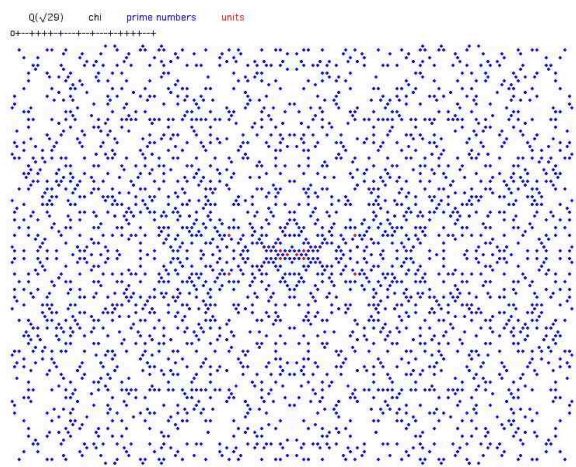
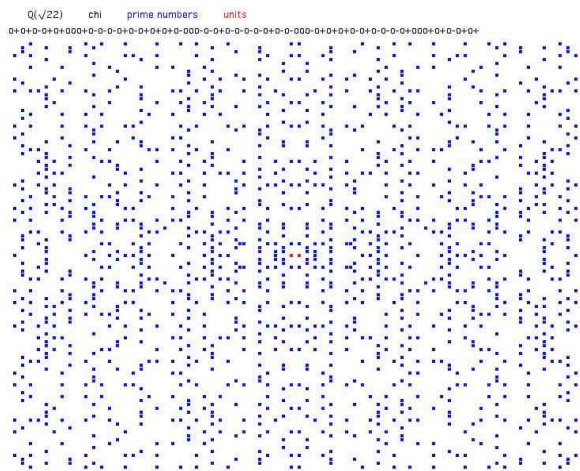
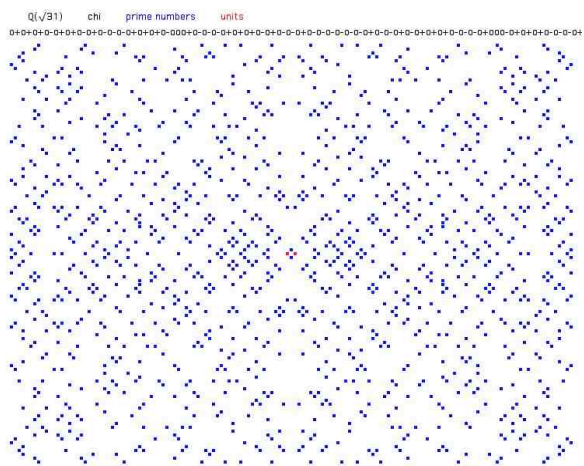
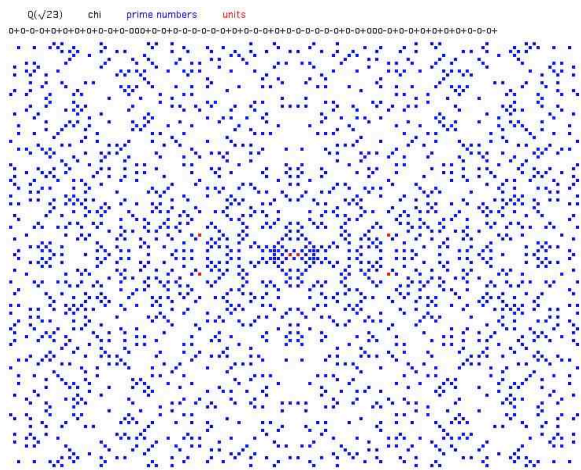

Appendix 3:	Pictures of prime numbers for complex non-UFD

The pictures show the quadratic character and a picture of prime numbers and units for some complex quadratic fields whose domain of integers is not a unique-factorization domain, namely of class numbers, h, as indicated

the fields of discriminant congruent 0 modulo 4:

h = 2: $Q(\sqrt{-5}), Q(\sqrt{-6}), Q(\sqrt{-10}), Q(\sqrt{-13})$, h = 4: $Q(\sqrt{-14}), Q(\sqrt{-17}), Q(\sqrt{-21})$,
h = 2: $Q(\sqrt{-22})$, h = 6: $Q(\sqrt{-26}), Q(\sqrt{-29})$

and the fields of discriminant congruent 1 modulo 4:

h = 2: $Q(\sqrt{-15})$, h = 3: $Q(\sqrt{-23}), Q(\sqrt{-31})$, h = 2: $Q(\sqrt{-35})$, h = 4: $Q(\sqrt{-39})$.

The pictures display the prime numbers, which generate the principal prime ideals, but not those irreducible numbers which are not prime. Moreover, the non-principal prime ideals are not displayed (see, however, appendices 5 and 7).

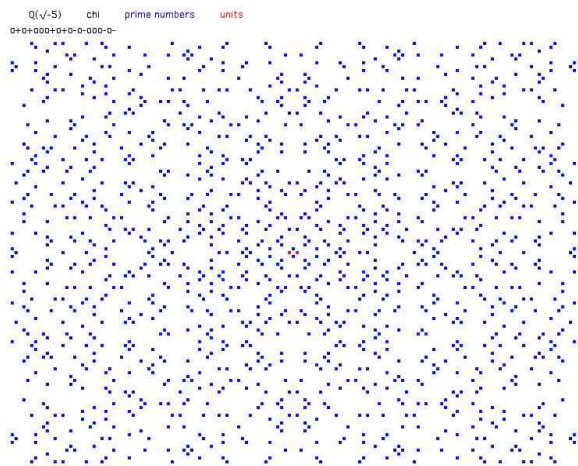
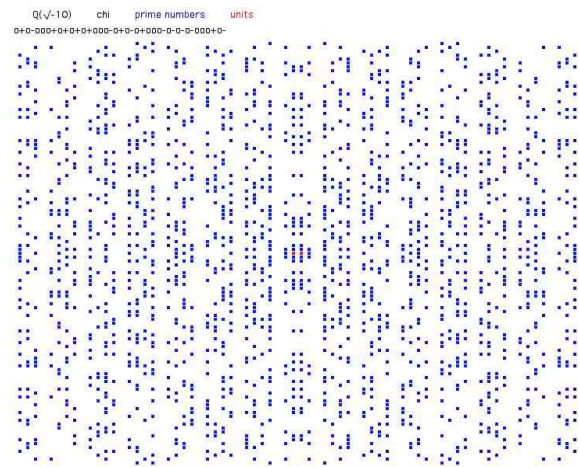
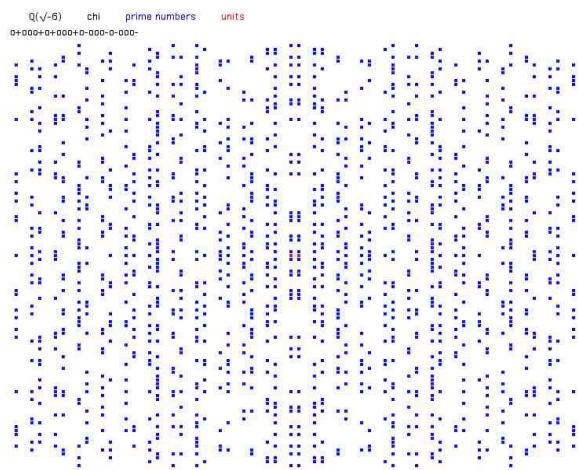
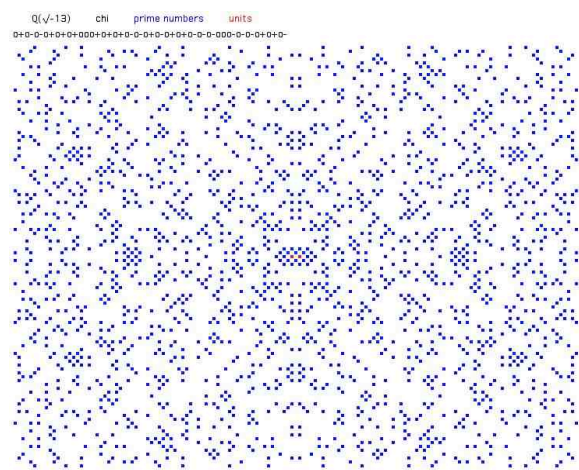

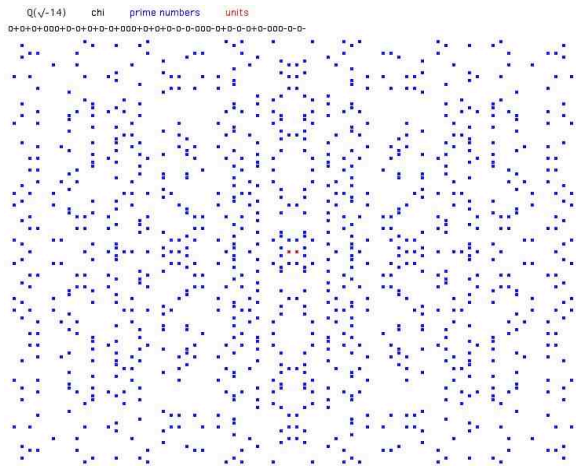
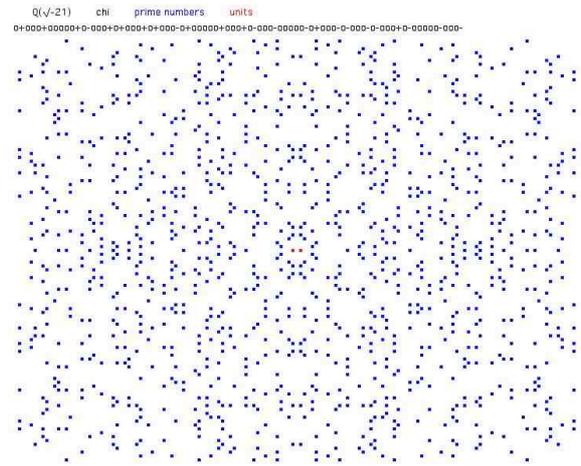
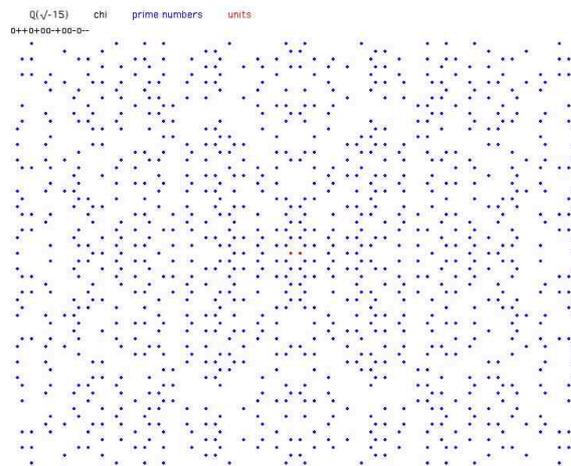
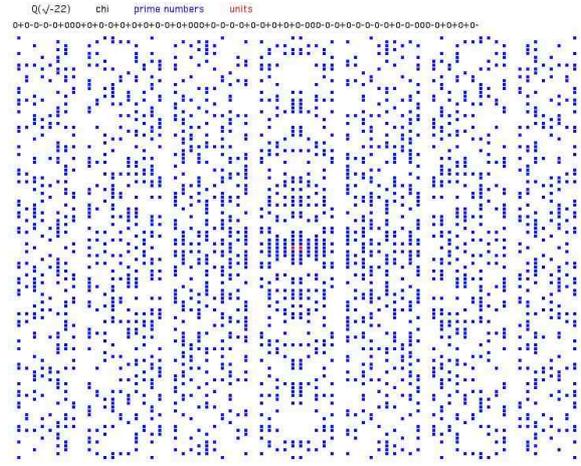
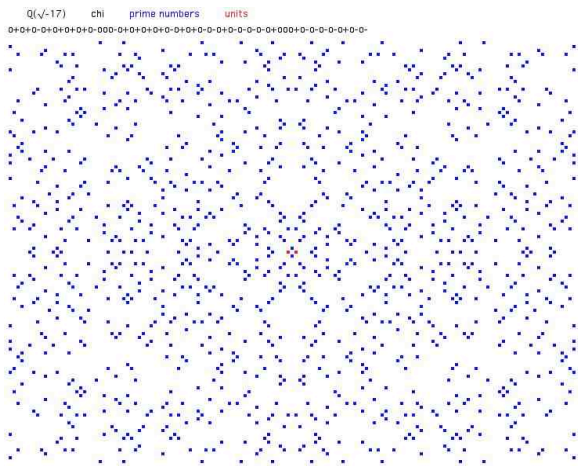
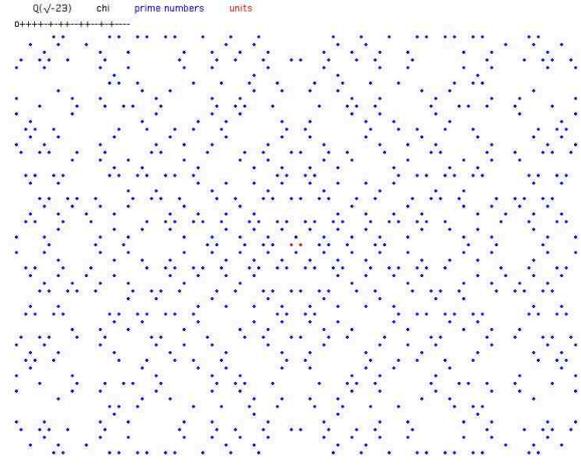

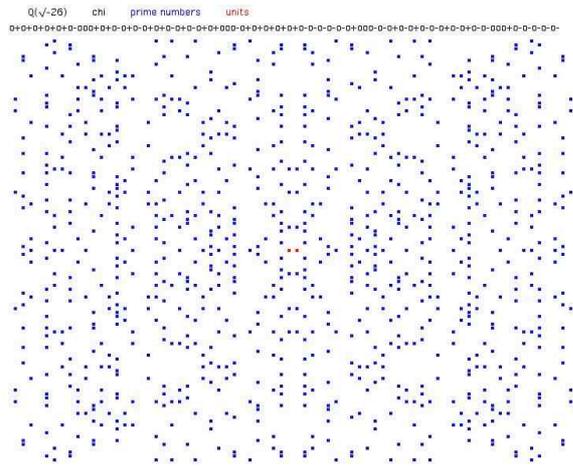
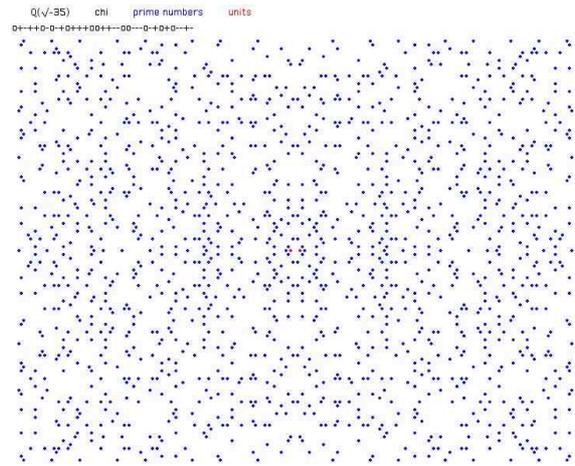
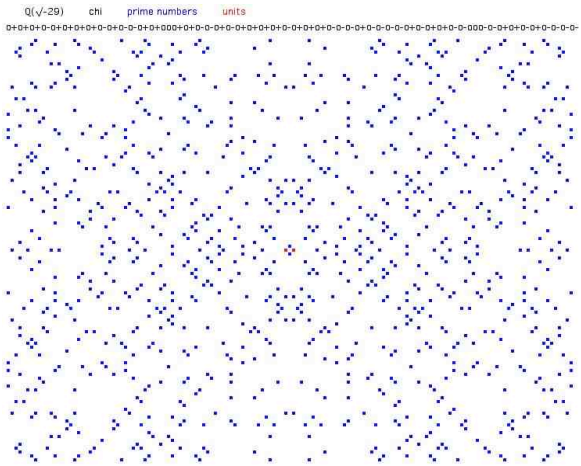
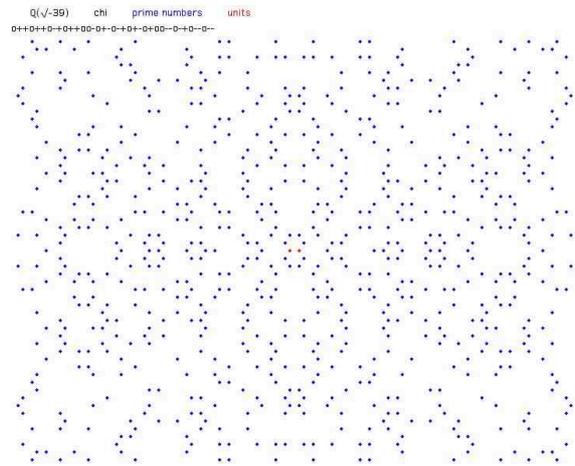
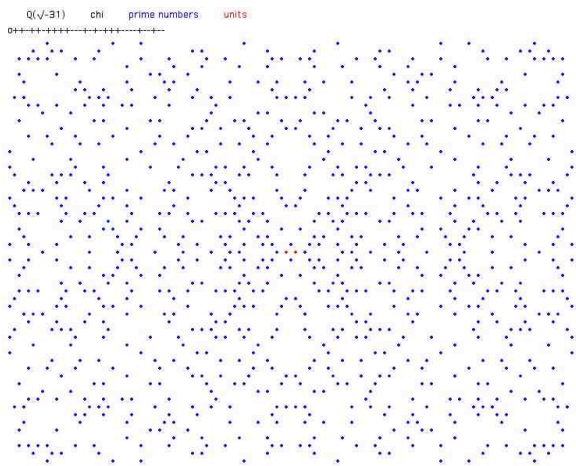

Appendix 4:    Pictures of prime numbers for real non-UFD

The pictures show the quadratic character and a picture of prime numbers and units for some complex quadratic fields whose domain of integers is not a unique-factorization domain, namely of class numbers, h, as indicated

the fields of discriminant congruent 0 modulo 4:

h = 2: Q($\sqrt{10}$), Q($\sqrt{15}$), Q($\sqrt{26}$), Q($\sqrt{30}$), Q($\sqrt{34}$), Q($\sqrt{35}$), Q($\sqrt{39}$), h = 3: Q($\sqrt{79}$), h = 4: Q($\sqrt{82}$)

and the fields of discriminant congruent 1 modulo 4:

h = 2: Q($\sqrt{65}$), Q($\sqrt{85}$), Q($\sqrt{105}$), h = 4: Q($\sqrt{145}$), h = 3: Q($\sqrt{229}$), Q($\sqrt{257}$).

The pictures display the prime numbers, which generate the principal prime ideals, but not those irreducible numbers which are not prime. Moreover, the non-principal prime ideals are not displayed (see, however, appendices 6 and 8).

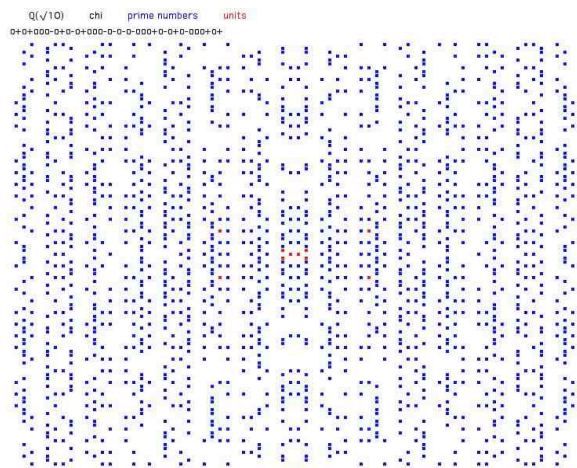
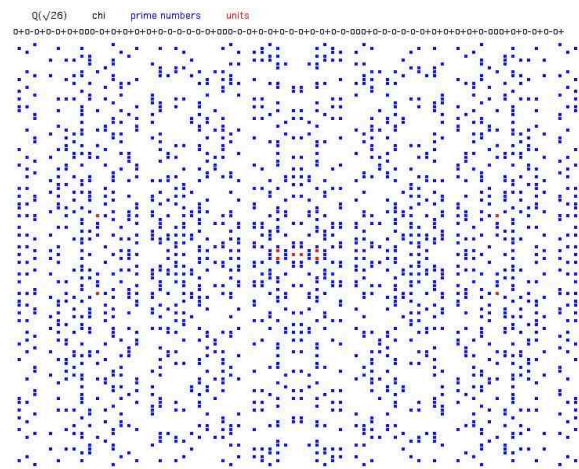
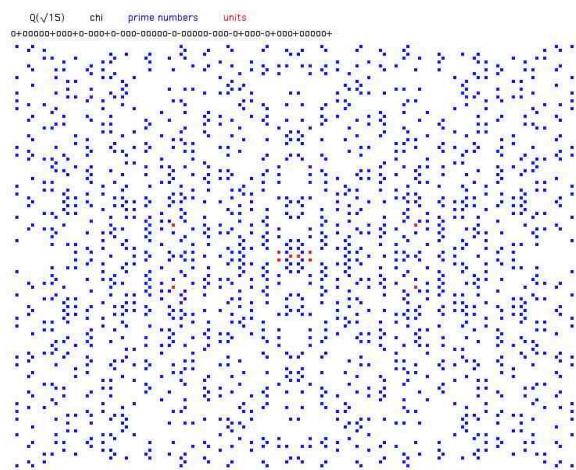
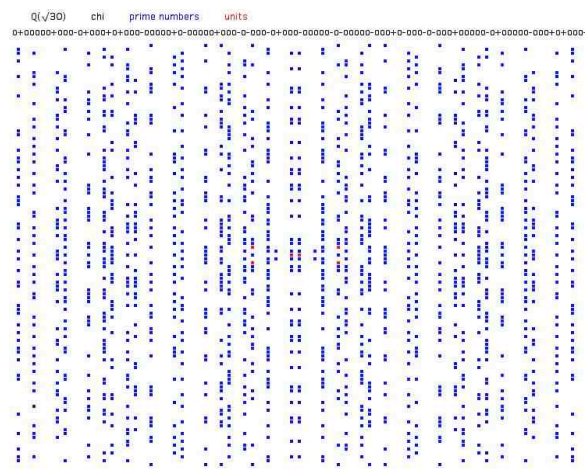

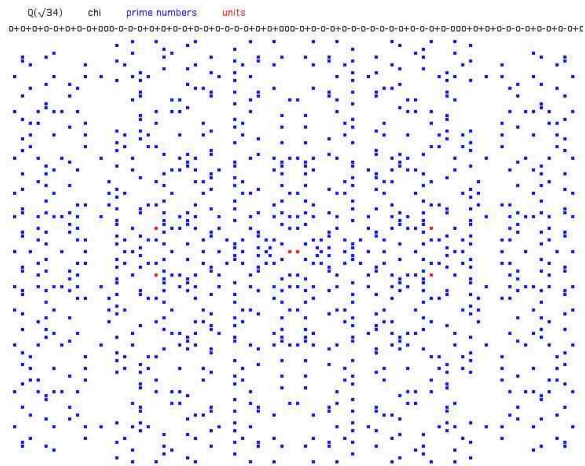
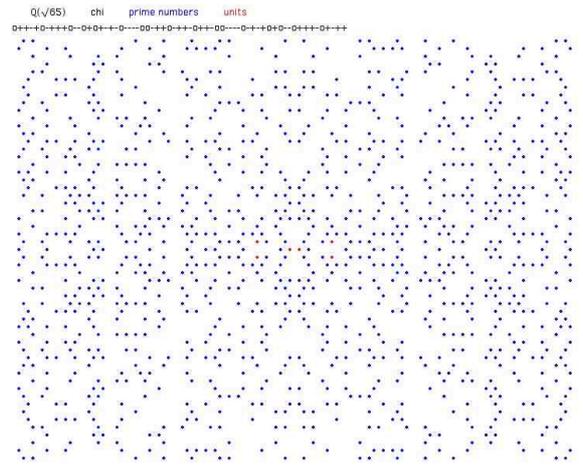
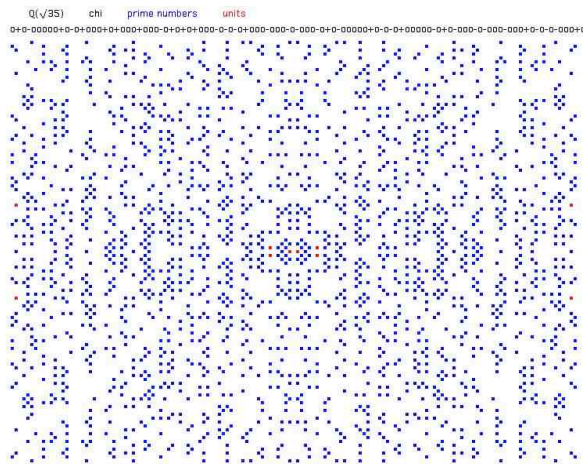
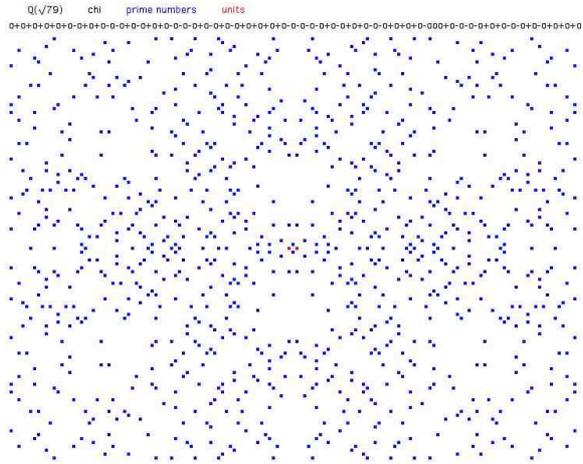
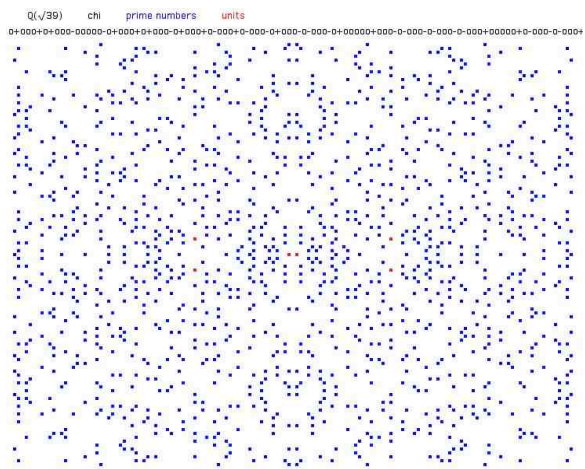
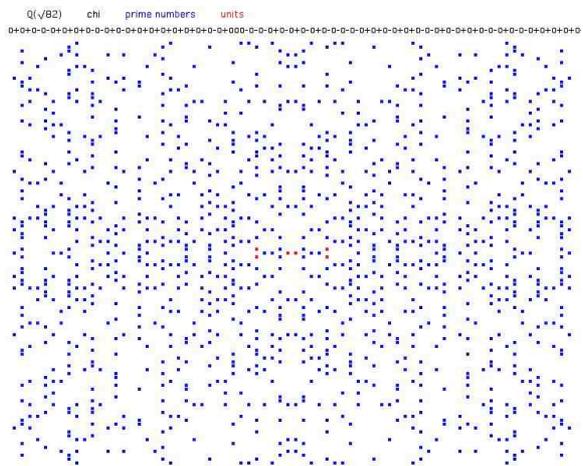

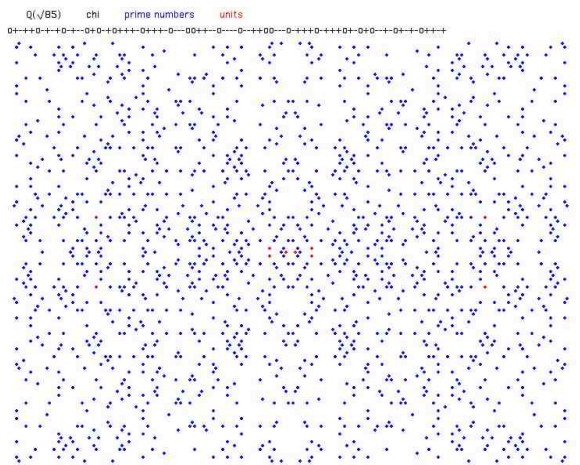
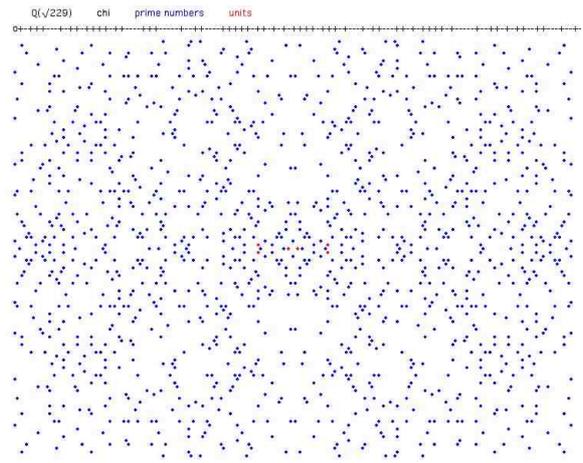
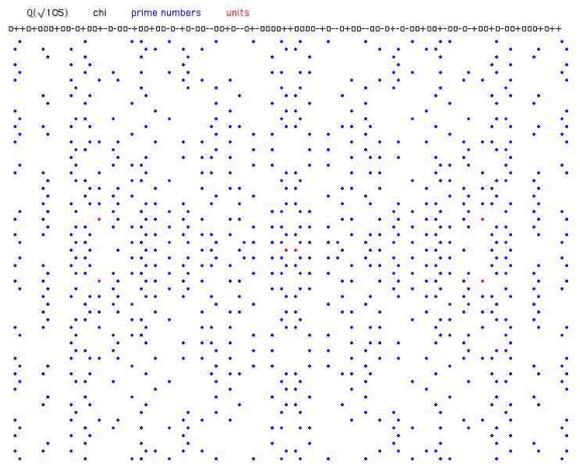
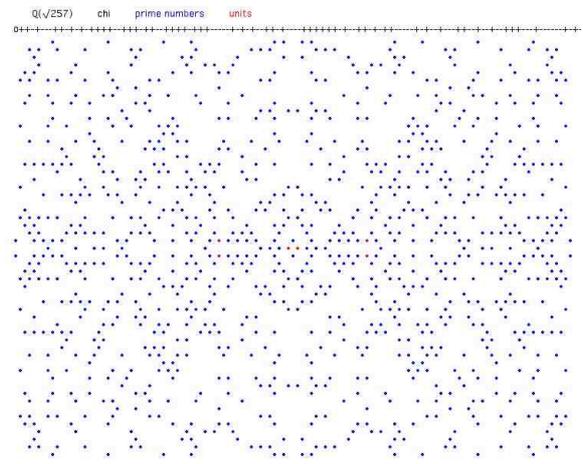
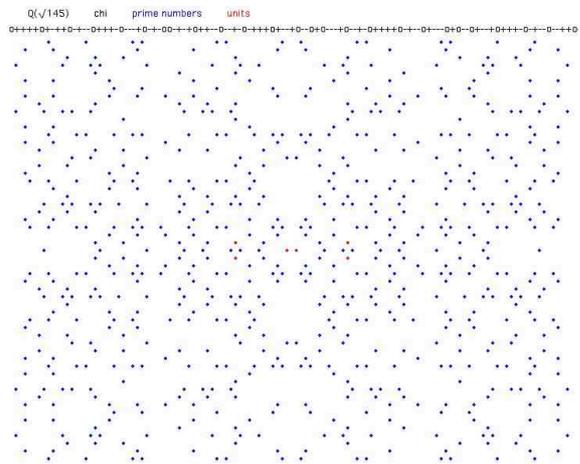

Appendix 5:    Pictures of prime numbers and ideals for complex fields of class number 2

The pictures show the quadratic character and a picture of prime numbers, units and non-principal prime ideals for some complex quadratic fields of class number 2, namely

the fields of discriminant congruent 0 modulo 4:

   Q(√-5), Q(√-6), Q(√-10), Q(√-13), Q(√-22)

and the fields of discriminant congruent 1 modulo 4:

   Q(√-15), Q(√-35), Q(√-51), Q(√-91), Q(√-115), Q(√-123), Q(√-187), Q(√-235).

The pictures display the prime numbers, which generate the principal prime ideals, but not those irreducible numbers which are not prime.
Moreover, the non-principal prime ideals are displayed as follows.
The non-principal ideals are obtained by dividing principal ideals by a certain non-principal prime ideal, I, generated by its norm and some integer of Q(√r). In the picture, the non-principal prime ideals then are represented by those numbers whose norm is equal to a prime norm times the norm of I. This norm of I is mentioned at the top of the picture.

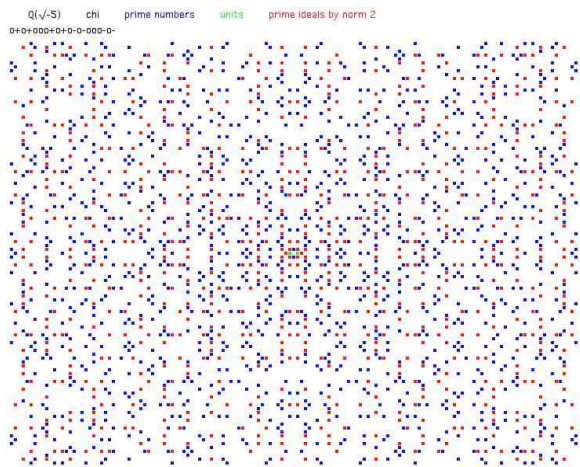
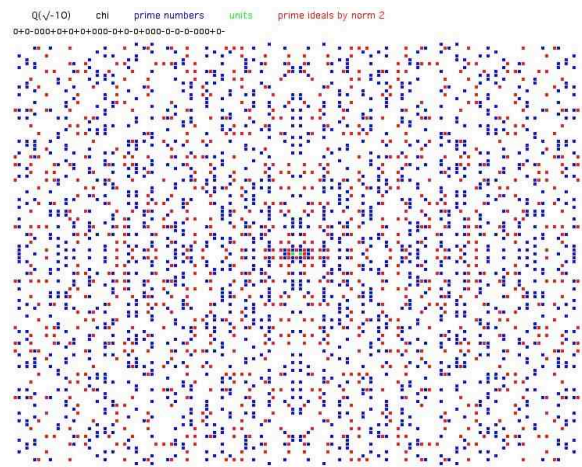
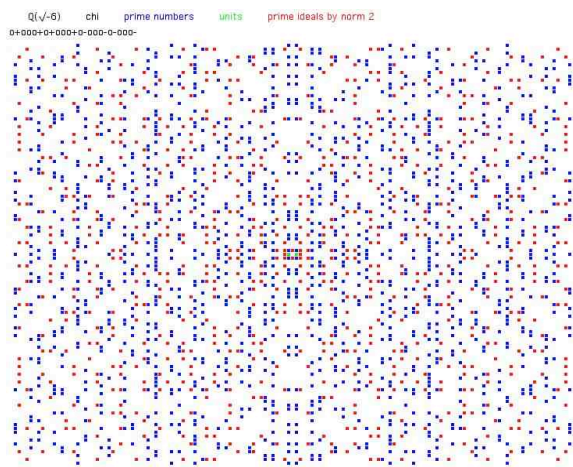
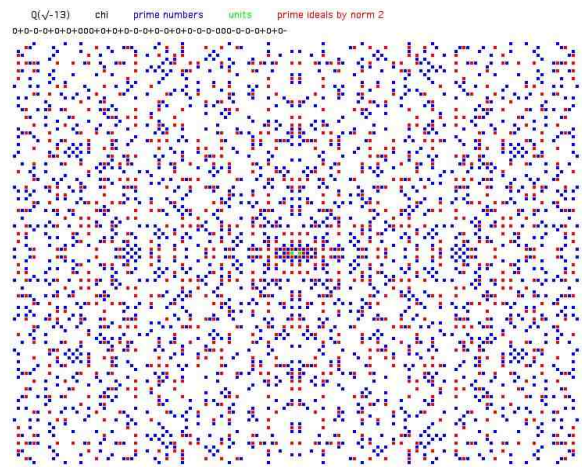

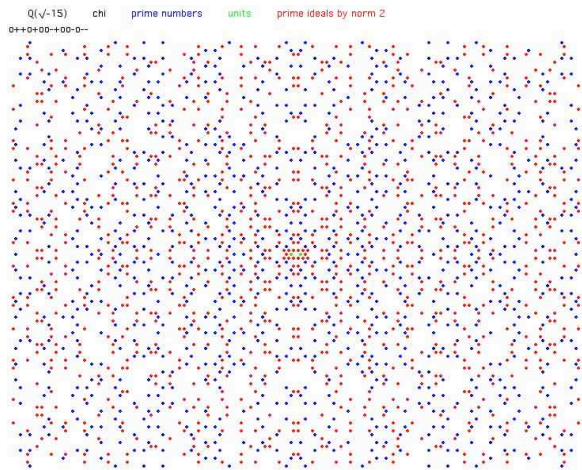
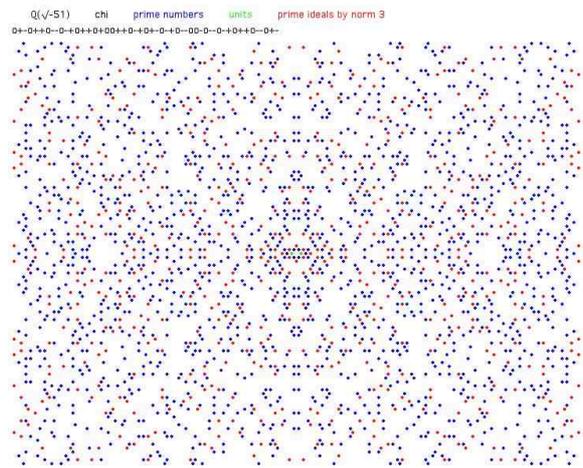
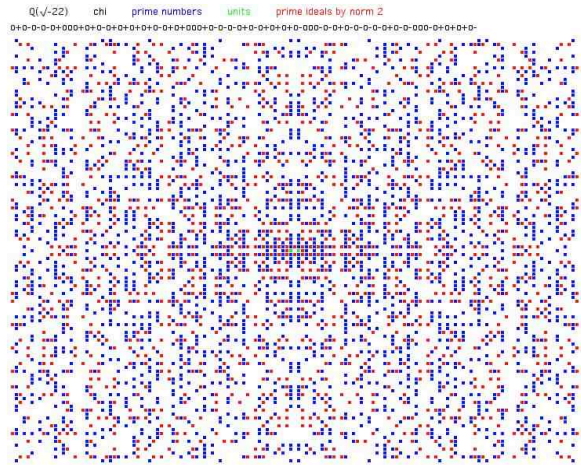
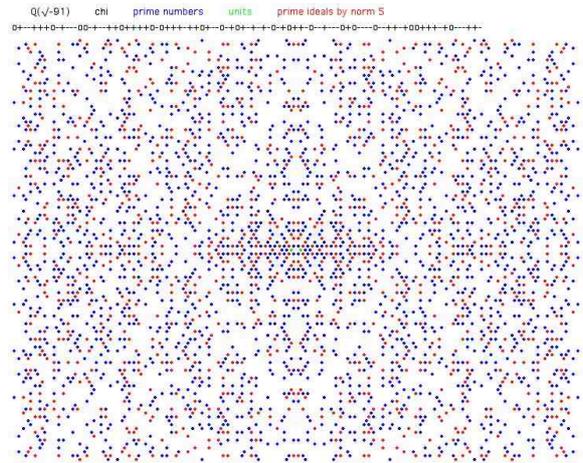
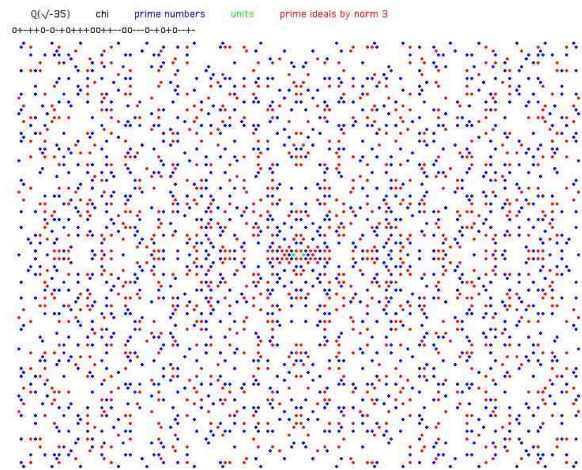
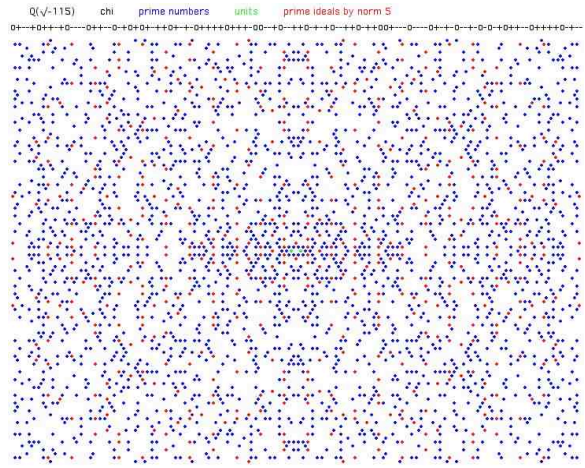

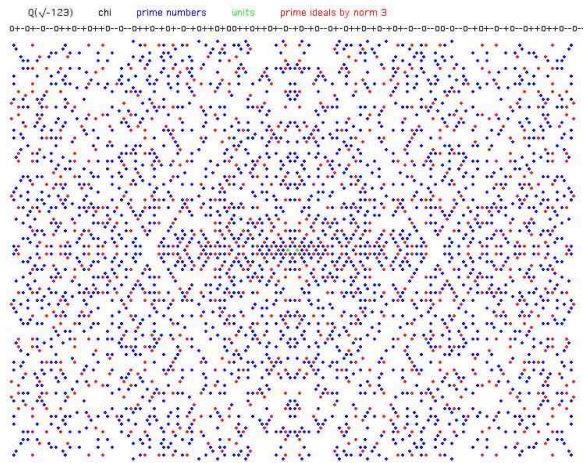
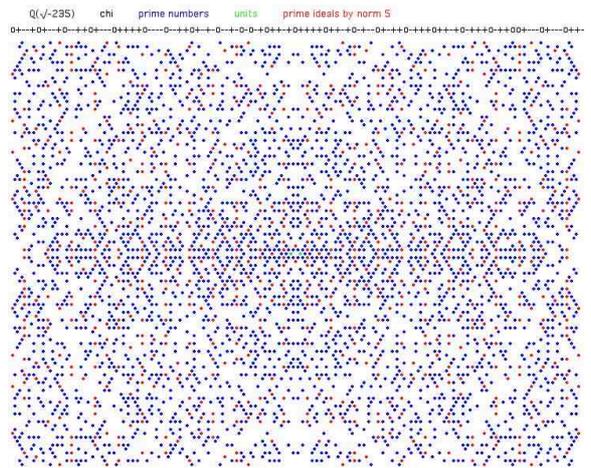
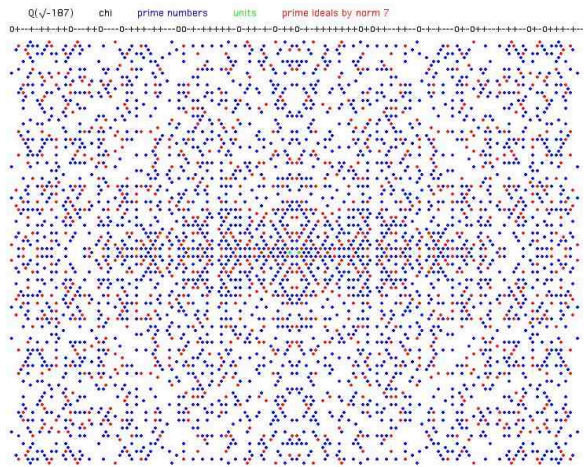

Appendix 6:    Pictures of prime numbers and ideals for real fields of class number 2

The pictures show the quadratic character and a picture of prime numbers, units and non-principal prime ideals for some real quadratic fields of class number 2, namely

the fields of discriminant congruent 0 modulo 4:

$Q(\sqrt{10}), Q(\sqrt{15}), Q(\sqrt{26}), Q(\sqrt{30}), Q(\sqrt{34}), Q(\sqrt{35}), Q(\sqrt{39})$

and the fields of discriminant congruent 1 modulo 4:

$Q(\sqrt{65}), Q(\sqrt{85}), Q(\sqrt{105})$.

The pictures display the prime numbers, which generate the principal prime ideals, but not those irreducible numbers which are not prime.
Moreover, the non-principal prime ideals are displayed as follows.
The non-principal ideals are obtained by dividing principal ideals by a certain non-principal prime ideal, I, generated by its norm and some integer of $Q(\sqrt{r})$. In the picture, the non-principal prime ideals then are represented by those numbers whose norm is equal to a prime norm times the norm of I. This norm of I is mentioned at the top of the picture.

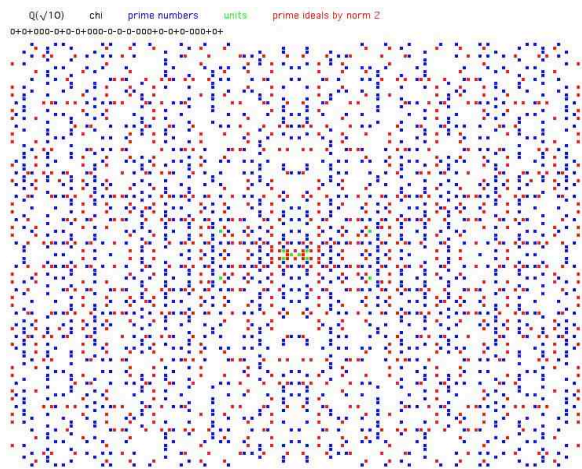
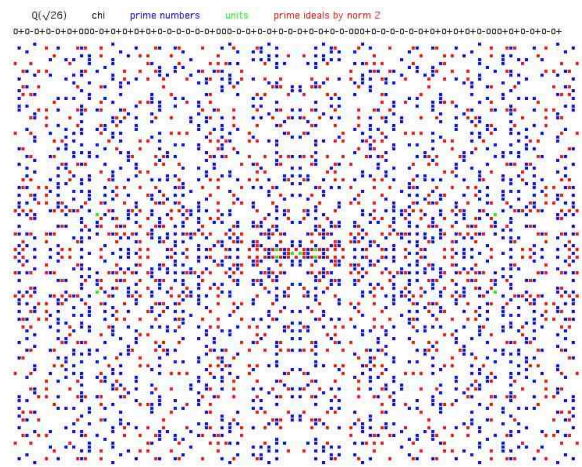
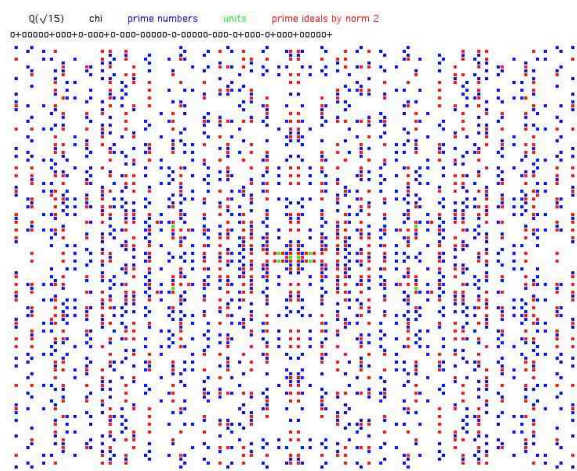
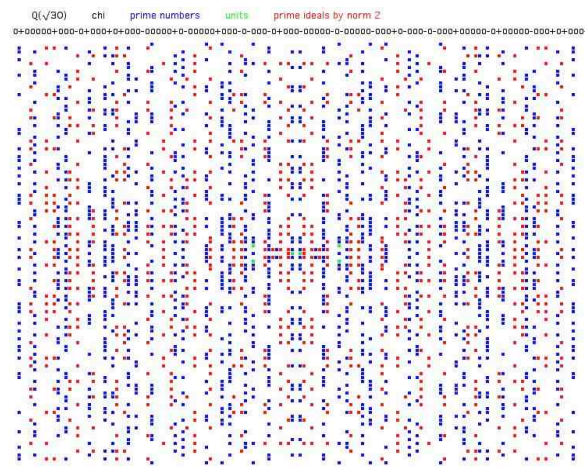

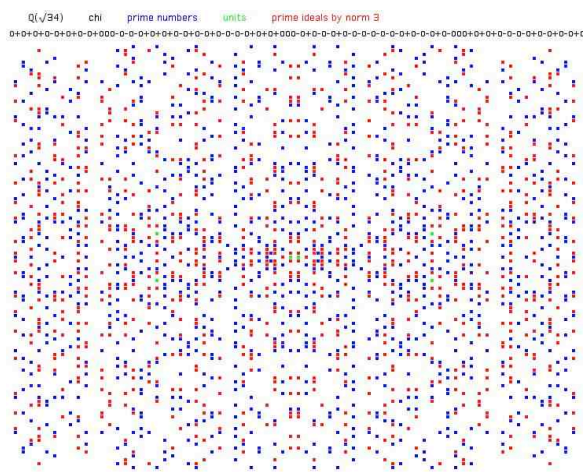
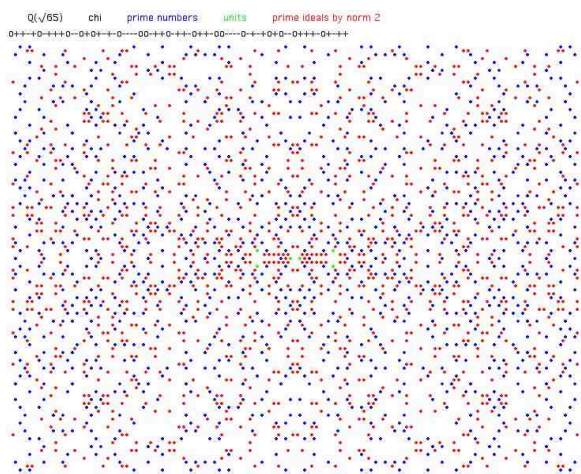
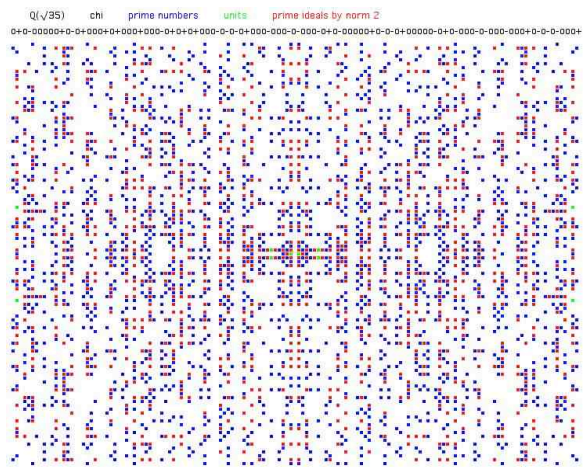
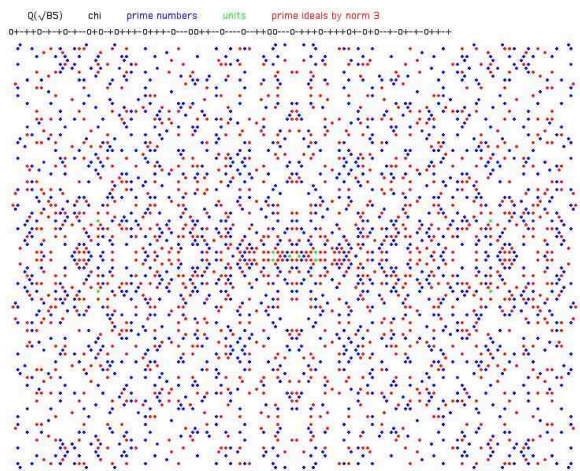
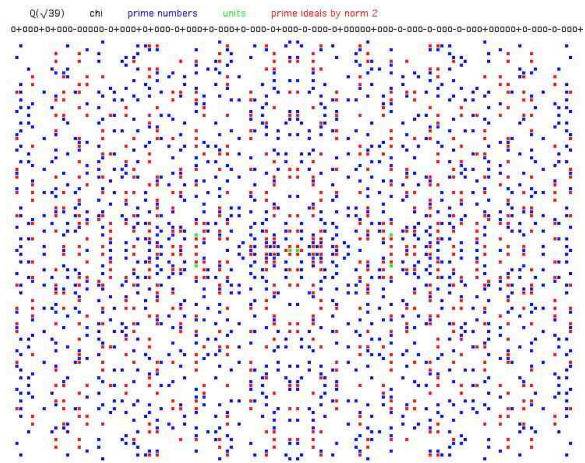
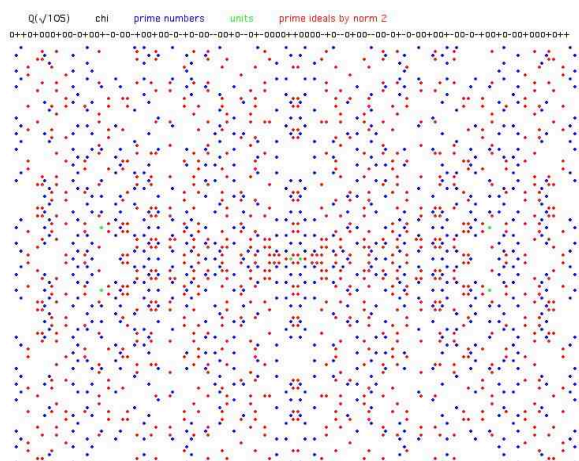

Appendix 7:　　Pictures of prime numbers and ideals for complex fields of class number 3

The pictures show the quadratic character and a picture of prime numbers, units and two mutually conjugate classes of non-principal prime ideals, one class red, and the other class green for the complex quadratic fields of class number 3 of discriminant larger than -350, namely
- since no complex fields of discriminant congruent 0 modulo 4 exist in this class -

the fields of discriminant congruent 1 modulo 4:

$$Q(\sqrt{-23}), Q(\sqrt{-31}), Q(\sqrt{-59}), Q(\sqrt{-83}), Q(\sqrt{-107}), Q(\sqrt{-139}), Q(\sqrt{-211}), Q(\sqrt{-283}),$$
$$Q(\sqrt{-307}), Q(\sqrt{-331}).$$

The pictures display the prime numbers, which generate the principal prime ideals, but not those irreducible numbers which are not prime.
Moreover, the non-principal prime ideals are displayed as follows.
The non-principal ideals are obtained by dividing principal ideals by a certain non-principal prime ideal, I, or its conjugate, where I := [norm, ζ], ζ := shift + (1 + √d) / 2,
i.e. I is generated by 'norm' being its norm, and the integer ζ of Q(√r).
In the picture, the non-principal prime ideals then are represented by those numbers whose norm is equal to a prime norm times the norm of I. This norm of I and shift are mentioned at the top of the picture, shift being needed to distinguish between the two mutually conjugate classes of non-principal ideals.

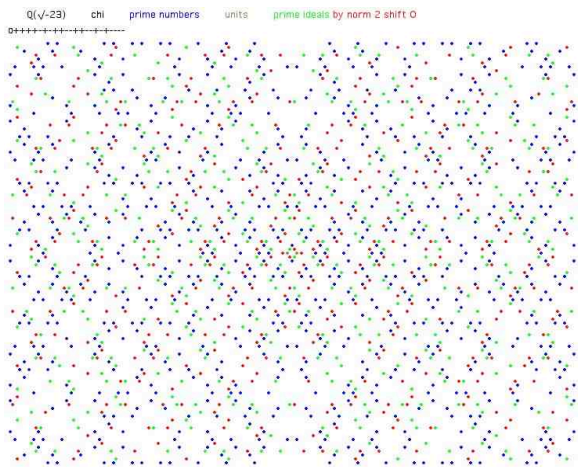
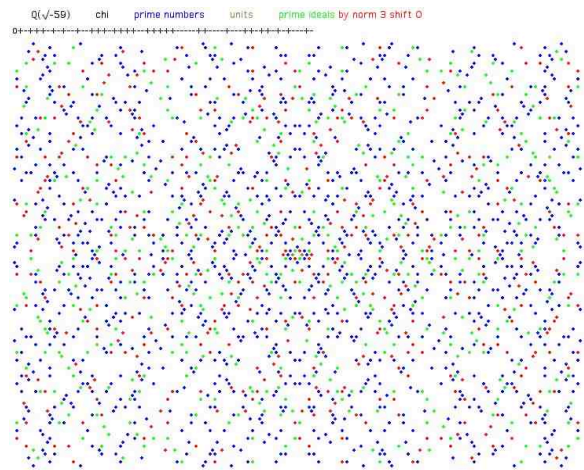
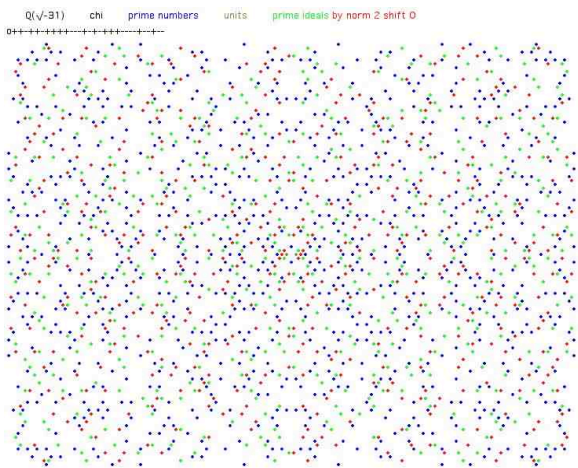
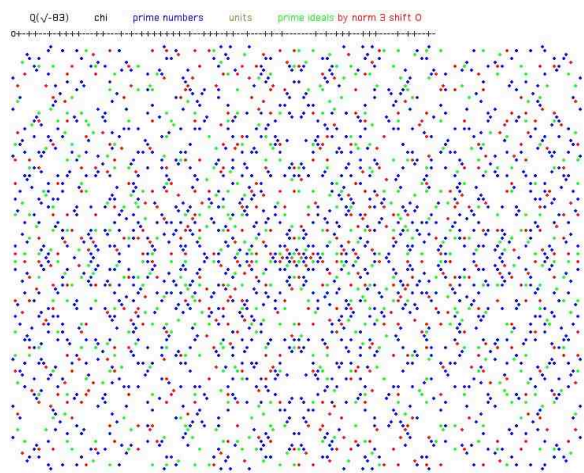

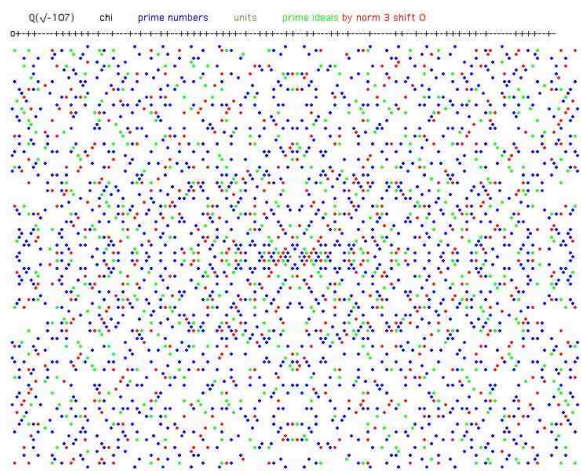
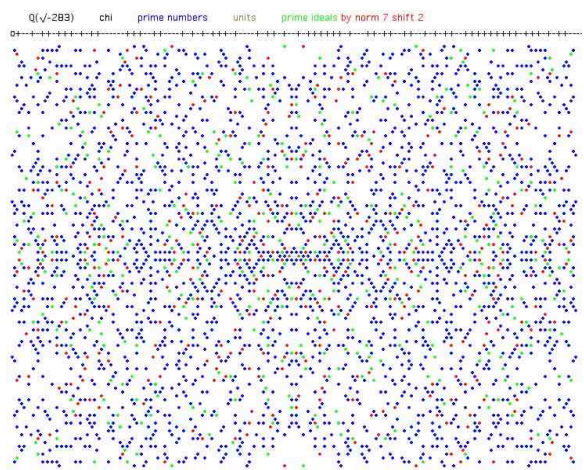
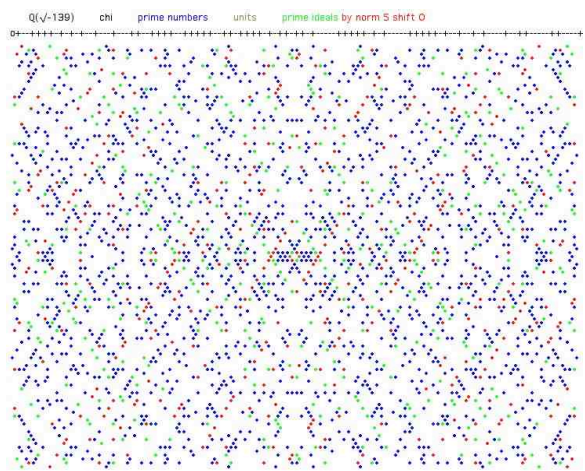
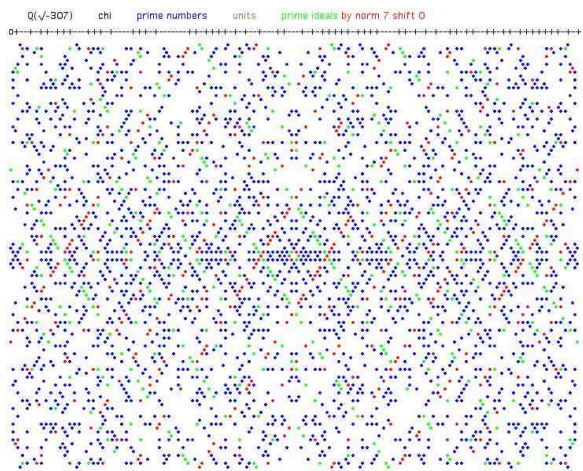
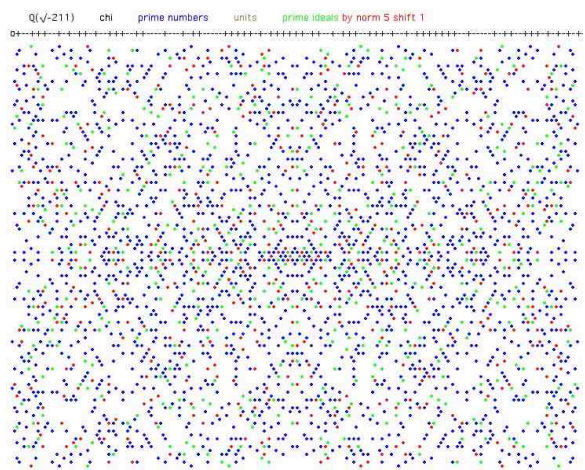
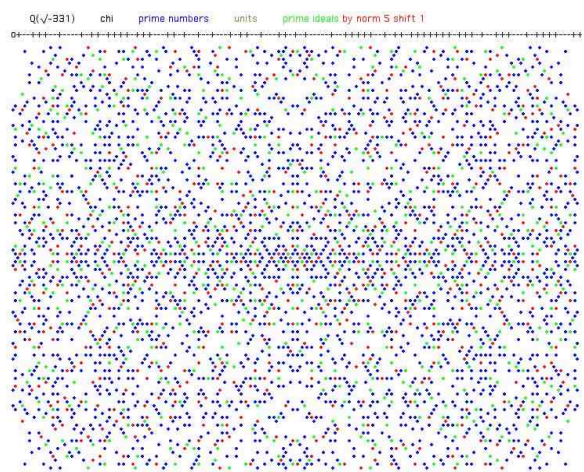

Appendix 8:            Pictures of prime numbers and ideals for real fields of class number 3

The pictures show the quadratic character and a picture of prime numbers, units and two mutually conjugate classes of non-principal prime ideals, one class red, and the other class green for some real quadratic fields of class number 3, namely

the fields of discriminant congruent 0 modulo 4:

   Q(√79), Q(√142), Q(√223), Q(√254), Q(√326), Q(√359)

and the fields of discriminant congruent 1 modulo 4:

   Q(√229), Q(√257), Q(√321).

The pictures display the prime numbers, which generate the principal prime ideals, but not those irreducible numbers which are not prime.
 Moreover, the non-principal prime ideals are displayed as follows.
The non-principal ideals are obtained by dividing principal ideals by a certain non-principal prime ideal, I, or its conjugate, where I := [norm, ζ], ζ := shift + (d mod 4 + √d) / 2,
i.e. I is generated by 'norm' being its norm, and the integer ζ of Q(√r).
In the picture, the non-principal prime ideals then are represented by those numbers whose norm is equal to a prime norm times the norm of I. This norm of I and shift are mentioned at the top of the picture, shift being needed to distinguish between the two mutually conjugate classes of non-principal ideals.

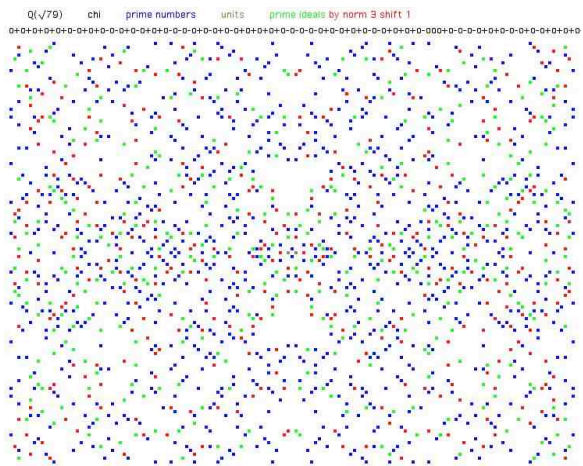
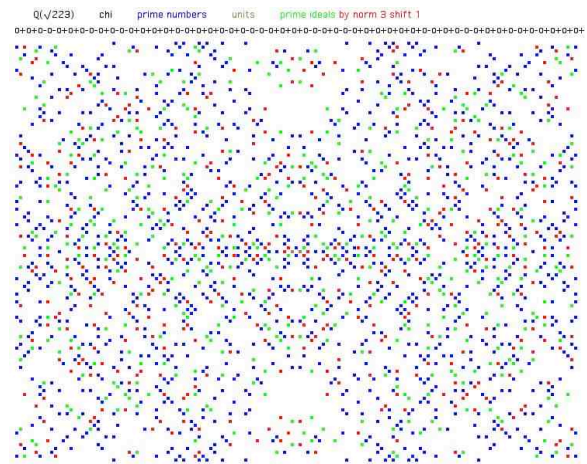
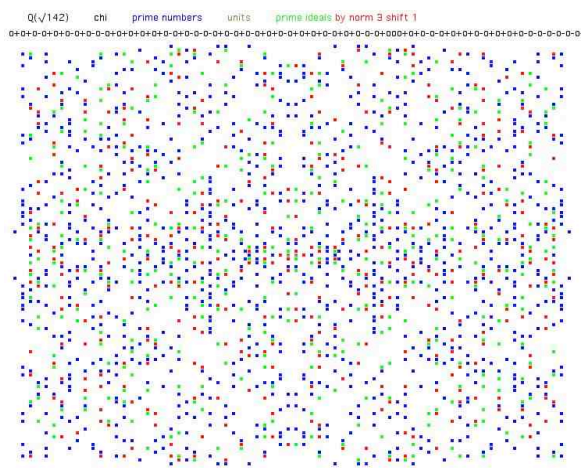
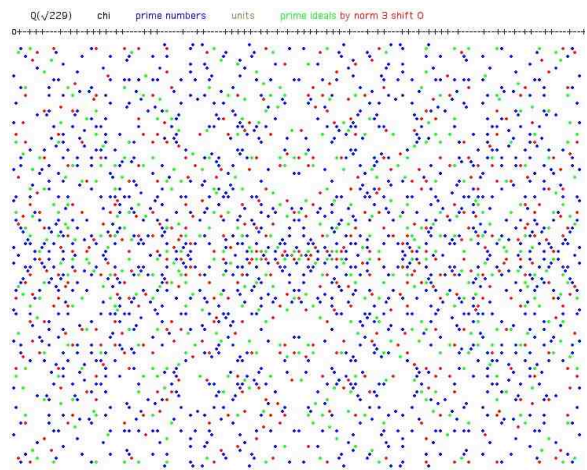

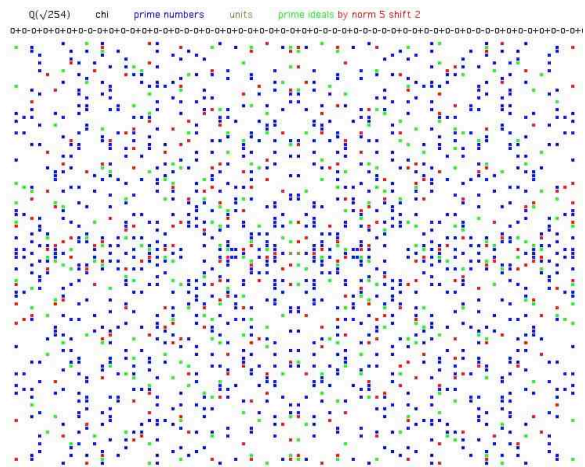
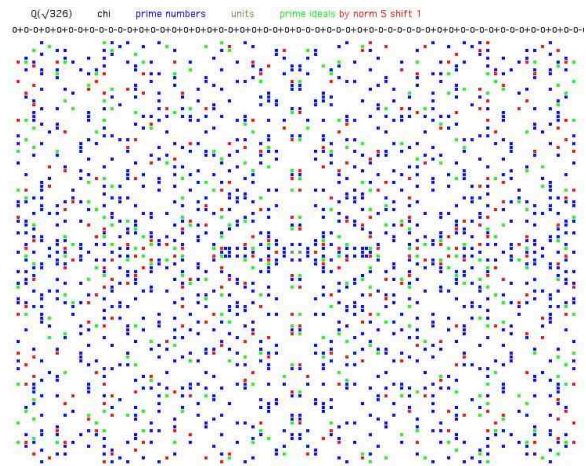
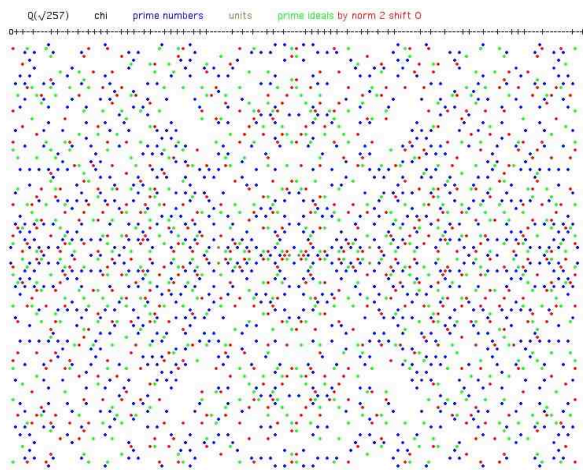
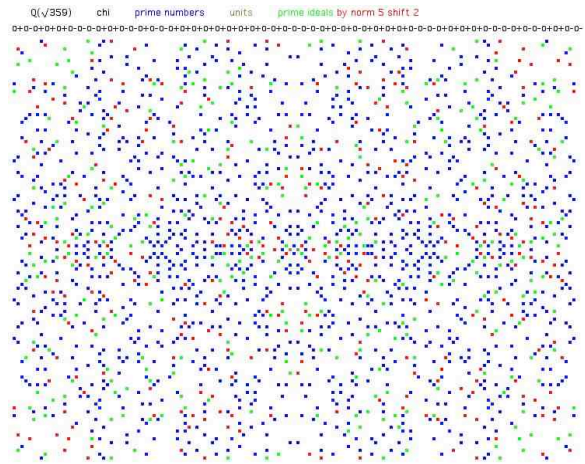
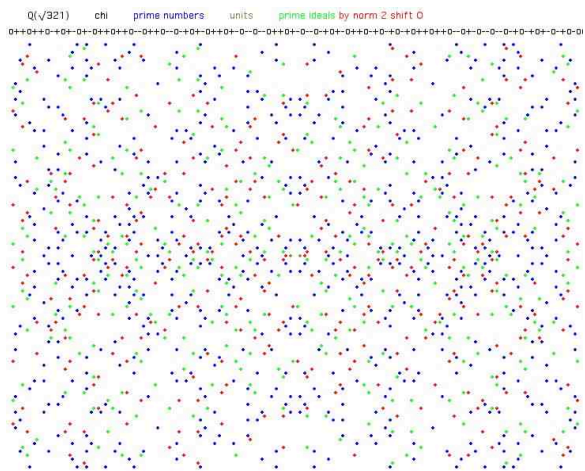